\documentclass[12pt]{article}
\usepackage{amssymb}
\newcommand{\sect}[1]{\setcounter{equation}{0}\section{#1}}

\newcommand{\be}{\begin{equation}}
\newcommand{\ee}{\end{equation}}
\newcommand{\bea}{\begin{eqnarray}}
\newcommand{\eea}{\end{eqnarray}}
\newcommand{\beano}{\begin{eqnarray*}}
\newcommand{\eeano}{\end{eqnarray*}}
\newcommand{\nonu}{\nonumber \\}

\newcommand{\hs}[1]{\hspace{#1 mm}}
\newcommand{\eps}{\epsilon}

\newcommand{\lbda}{\lambda}

\newcommand{\ca}{\mbox{$\cal{A}$}}

\newcommand{{\cg}}{\mbox{$\cal{G}$}}

\newcommand{\cm}{\mbox{$\cal{M}$}}

\newcommand{\ct}{\mbox{$\cal{T}$}}
\newcommand{\cu}{\mbox{${\cal U}$}}
\newcommand{\cw}{\mbox{$\cal{W}$}}
\newcommand{\cy}{\mbox{${\cal Y}$}}

\newcommand{\wt}[1]{\widetilde{#1}}
\newcommand{\mb}[1]{\hs{4}\mbox{#1}\hs{4}}
\newcommand{\half}{\frac{1}{2}}

\newtheorem{theor}{Theorem}[section]
\newtheorem{prop}[theor]{Property}
\newtheorem{defi}[theor]{Definition}

\newtheorem{coro}[theor]{Corollary}
\newcommand{\prf}{\underline{Proof:}\ }
\newcommand{\finprf}{\null \hfill {\rule{5pt}{5pt}}\\[2.1ex]\indent}
\newcommand{\ie}{{\it i.e.}\ }

\newtheorem{guess}{Conjecture}
\newcommand{\CC}{\mbox{${\mathbb C}$}}

\newcommand{\ZZ}{\mbox{${\mathbb Z}$}}
\newcommand{\NN}{\mbox{${\mathbb N}$}}

\newcommand{\II}{\mbox{${\mathbb I}$}}

\newcommand{\bara}{{\bar{a}}}
\newcommand{\barb}{{\bar{b}}}
\newcommand{\barc}{{\bar{c}}}
\newcommand{\bard}{{\bar{d}}}
\newcommand{\bare}{{\bar{e}}}
\newcommand{\bari}{{\bar{\imath}}}
\newcommand{\barj}{{\bar{\jmath}}}

\newcommand{\tilm}{{\widetilde{m}}}

\newcommand{\dg}{\mbox{deg}}

\newcommand{\EE}{\mbox{${\mathbb E}$}}
\newcommand{\FF}{\mbox{${\mathbb F}$}}

\newcommand{\JPhys}[1]{Journ.\ Phys.\ {\bf #1}}
\newcommand{\LMP}[1]{Lett.\ Math.\ Phys.\ {\bf #1}}
\newcommand{\IJMP}[1]{Int. J.\ Mod.\ Phys.\ {\bf #1}}

\begin{document}
\renewcommand{\thefootnote}{\fnsymbol{footnote}}
\newpage
\pagestyle{empty}
\setcounter{page}{0}

\newcommand{\LAP}{LAPTH}
\def\logo{{\bf {\huge LAPTH}}}

\centerline{\logo}

\vspace {.3cm}

\centerline{{\bf{\it\Large 
Laboratoire d'Annecy-le-Vieux de Physique Th\'eorique}}}

\centerline{\rule{12cm}{.42mm}}

\vspace{20mm}

\begin{center}

  {\LARGE  {\sffamily Twisted superYangians and their representations
    }}\\[1cm]

\vspace{10mm}
  
{\large C. Briot and E. Ragoucy\footnote{ragoucy@lapp.in2p3.fr}\\[.42cm]
   Laboratoire de Physique Th{\'e}orique \LAP\footnote{URA 14-36 
    du CNRS, associ{\'e}e {\`a} l'Universit{\'e} de Savoie.}\\[.242cm]
    LAPP, BP 110, F-74941  Annecy-le-Vieux Cedex, France. }
\end{center}
\vfill

\begin{abstract}
Starting with the superYangian $Y(M|N)$ based on $gl(M|N)$, 
we define twisted superYangians $Y(M|N)^\pm$. Only $Y^+(M|2n)$ and 
$Y^-(2m|N)$ can be defined, and appear to be 
isomorphic one with each other. We study their
finite-dimensional irreducible highest weight representations.
\end{abstract}
\vfill
\vfill

\rightline{\tt math.QA/yymmnn}
\rightline{\LAP-875/01}
\rightline{November 2001}

\newpage
\pagestyle{plain}
\setcounter{footnote}{0}

\markright{\today\dotfill DRAFT\dotfill }


\sect{Introduction}
Quite recently, a revival of interest has been put on coideal algebras of
Hopf algebras, both from mathematical and physical point of view. Among these 
algebras, let us note the twisted Yangians $Y^\pm(N)$, introduced by Olshanski
 \cite{olsh} and widely studied (see for instance \cite{rapM} and references 
therein), or the reflection algebras, introduced by Sklyanin \cite{skly} 
and studied in \cite{BNLS,MR}.

From the 
mathematical point of view, it seems that such coideal condition is quite 
restrictive, leading to a very small class of subalgebras, for a given Hopf 
algebra \cite{letz,noumi}. Indeed, for quantum algebras $\cu_q(gl_N)$, it has 
been proven that they are natural deformations of symmetric spaces 
\cite{letz}.

From a physical point of view, such ideals appear to play an important 
role in integrable systems with boundaries \cite{BNLS, delius, OV-bound}.
They appear to be the integrals of motion of such systems 
\cite{BNLS, OV-bound}, and also naturally deduced from the boundary 
condition \cite{delius}.

\null

It appears thus natural to look for such coideals when the 
underlying algebra is not $gl_N$ anymore. Such types of algebras have 
been introduced in 
\cite{AACFR} for the case of Yangians based on orthogonal and symplectic 
algebras, and orthosymplectic superalgebras. They are defined as the 
``twist'' of the (super)Yangian based on the corresponding Lie 
(super)algebra. 

The aim of the present work is to complete the picture 
with the case of superYangians based on $gl(M|N)$.
After recalling the basic definitions and properties of the superYangians 
$Y(gl(M|N))\equiv Y(M|N)$ in section \ref{supY}, we will define the 
twisted superYangians $Y(M|N)^+$ in section \ref{superYtw}. Their 
finite-dimensional irreducible representations are studied in 
section \ref{irreps}. We conclude in 
section \ref{conc}.

\sect{Super Yangians $Y(M|N)$\label{supY}}
The superYangian $Y(M|N)$ based on the $gl(M|N)$ superalgebra
has been introduced in \cite{naz}, and its irreducible finite-dimensional 
representations studied in 
\cite{zhang}. Since it is a $\ZZ_{2}$-graded (Hopf) algebra, different 
conventions can be chosen: the ones we choose are given below. We will 
rephrase the properties given in \cite{zhang} in this context. 

\subsection{Graded spaces\label{graded}}
We start with $K\times K$ matrices acting on the vector space 
$\CC^K$, and introduce a $\ZZ_{2}$-grading $[.]$ on these spaces.
We will denote by $E_{ij}$ the usual $K\times K$ matrices which have 
1 in position $(i,j)$, and $e_{i}$ the basic vectors of $\CC^{K}$ 
which have 1 in position $i$.
\be
E_{ij}e_{k}=\delta_{jk}e_{i}
\ee
The $\ZZ_{2}$-grade is defined by
\be
[E_{ij}]=[i]+[j] \ \mbox{ ; }\ [e_{i}]=[i]\ \mbox{ and }\ 
[i]\in\{0,1\}\ \forall i,j=1,\ldots,K
\ee
We will called even the matrices and vectors such that
\be
A=A_{ij}E_{ij}\mb{with} [A_{ij}]=[i]+[j]
\mb{;} u=u_{i}e_{i}\mb{with} [u_{i}]=[i]
\ee

The tensor product of graded matrices is chosen to be graded:
\be
(E_{ij}\otimes E_{kl})(E_{ab}\otimes E_{cd}) = (-1)^{([k]+[l])([a]+[b])}
(E_{ij}E_{ab})\otimes (E_{kl}E_{cd})
\ee
On tensor product of $\CC^{K}$ vector spaces, one has
\be
(E_{ij}\otimes E_{kl})(e_a\otimes e_b)=(-1)^{([k]+[l])[a]\,[u]}\,
(E_{ij}e_a)\otimes(E_{kl}e_b)
\ee
We introduce the graded permutation operator:
\be
P_{12}=\sum_{{i,j}}(-1)^{[j]} E_{ij}\otimes E_{ji}
\ee
which obeys $P^{2 }=\II$ and 
\be
P (e_i\otimes e_j) =  (-1)^{[i][j]}e_j\otimes e_i \mbox{ , }
P (E_{ij}\otimes E_{kl})P = (-1)^{([i]+[j])([k]+[l])}E_{kl}\otimes E_{ij}
\ee
\subsection{Definition and first properties of $Y(M|N)$}
We set $K=M+N$ and define the $\ZZ_{2}$-grade by:
\be
\begin{array}{l}
{[i]}=0\ \mbox{ for }\ 1\leq i\leq M\\
{[i]}=1\ \mbox{ for }\ M+1\leq i\leq M+N
\end{array}
\ee
The super Yangian $Y(M|N)$ has generators $T^{ab}_{(n)}$ (of 
$\ZZ_{2}$-grade $[a]+[b]$), gathered in:
\be
T(u)=\sum_{a,b=1}^K\sum_{n\geq0}u^{-n}\,T^{ab}_{(n)}E_{ab}=
\sum_{a,b=1}^KT^{ab}(u)E_{ab}=\sum_{n\geq0}u^{-n}\,T_{(n)}
\mb{with} T_{(0)}=\II_{K}
\ee
The even matrix $T(u)\in M_{K}[Y(M|N)]$ obeys:
\bea
&&R_{12}(u-v)\, T_{1}(u)\, T_{2}(v) = T_{2}(v)\, T_{1}(u)\, R_{12}(u-v)
\label{RTT}\\
&&\mbox{ with } R_{12}(x)=\II\otimes\II -\frac{1}{u-v}P_{12}\mb{,} 
P_{12}= \sum_{{i,j}}(-1)^{[j]} E_{ij}\otimes E_{ji}
\label{matR}
\eea
or equivalently
\be
{[T^{ab}(u)\,,\, T^{cd}(v)\}} = \frac{(-1)^{[a][b]+([a]+[b])[c]}}{u-v} 
\left(\rule{0ex}{2.4ex}T^{cb}(u) T^{ad}(v)-T^{cb}(v) T^{ad}(u)\right)
\label{comYMN}
\ee
where the graded commutator is defined by
\be
{[A\,,\, B\}}= A\, B -(-1)^{[A][B]}\, B\, A
\ee
It is a Hopf algebra:
\be
\Delta(T^{ab}(u))=\sum_{c=1}^{M+N}T^{ac}(u)\otimes T^{cb}(u)\ ;\ 
\eps(T^{ab}(u))=\delta^{ab}\ ;\ 
S(T^{ab}(u))=\left(T^{-1}(u)\right)^{ab} \label{DeltaT}
\ee
Note that $Y(M|N)$ contains several subalgebras\cite{zhang}:
\begin{prop}[Subalgebras of $Y(M|N)$]\hfill
    \label{subalg}
\null\ \\
    
$(i)$ The generators $T^{ab}(u)$ with $a,b=1,\ldots,M$ (resp. $T^{ab}(u)$, 
    $a,b=M+1,\ldots,M+N$) define the algebra inclusion $Y(M)\subset Y(M|N)$ 
    (resp. $Y(N)\subset Y(M|N)$).
    
$(ii)$ The generators $T^{ab}_{(1)}$ with $a,b=1,\ldots,M+N$ form a $gl(M|N)$ 
   Lie sub-super\-algebra of $Y(M|N)$.
    
$(iii)$ The generators $T^{ab}(u)$ with $a,b\in\{M,M+1\}$ define the algebra inclusion 
     $Y(1|1)\subset Y(M|N)$.
\end{prop}
{\bf Remark:}
The above subalgebras are \underline{not} Hopf subalgebras of $Y(M|N)$. 
Indeed, the coproduct (\ref{DeltaT}) is based on all the generators of 
$Y(M|N)$, so that it does not induce the coproduct of $Y(M)$ and 
$Y(N)$.

We have also the property
\begin{prop}[Isomorphism between $Y(M|N)$ and $Y(N|M)$]\hfill\\
\label{YMN-YNM}
Let
\be
\begin{array}{l}
\wt{T}^{ab}(u)=(-1)^{[\bara]([\barb]+1)} T^{\barb\bara}(u)
\mbox{ with }\bara=K+1-a\\[1.2ex]
[a]'=[\bara]+1
\end{array}
\ee
$\wt{T}(u)$ obey the Hopf algebra relations of $Y(N|M)$, with 
$\Delta'(x)=P\Delta(x)P$. We have thus 
a Hopf algebra isomorphism between $Y(M|N)$ and $Y(N|M)$.
\end{prop}
\prf
One first proves that $\wt{T}(u)$ satisfies the commutation relations 
of $Y(N|M)$. 

To prove that, we need:
\[
    T^{cb}(u) T^{ad}(v)-T^{cb}(v) T^{ad}(u)=- (-1)^{([a]+[d])([b]+[c])}
    \left(\rule{0ex}{2.1ex}T^{ad}(u) T^{cb}(v)-T^{ad}(v) T^{cb}(u)\right)
\]
which can be proven either by a direct calculation, or by using the 
graded antisymmetry of the commutator:
\be
{[T^{ab}(u)\,,\, T^{cd}(v)\}}=-(-1)^{([a]+[d])([b]+[c])}\, 
[T^{cd}(v)\,,\, T^{ab}(u)\}
\ee
and computing $[T^{cd}(v)\,,\, T^{ab}(u)\}$ using (\ref{comYMN}) with 
the replacements $(a,b,u)\leftrightarrow(c,d,v)$.

With the help of the above calculation
one gets:
\beano
{[\wt{T}^{ab}(u)\,,\, \wt{T}^{cd}(v)\}} &=& 
(-1)^{[\bara][\barb]+[\barc][\bard]}\
{[T^{\barb\bara}(u)\,,\, T^{\bard\barc}(v)\}}\\
&=& \frac{(-1)^{[\bara][\barb]+[\barc][\bard]+[\barb][\bara]+([\barb]+[\bara])[\bard]}}{u-v}\
\left(\rule{0ex}{2.1ex}T^{\bard\bara}(u) T^{\barb\barc}(v)-
T^{\bard\bara}(v) T^{\barb\barc}(u)\right)\\
&=& \frac{(-1)^{([\bara]+[\barb]+[\barc])[\bard]
+([\barb]+[\barc])([\bara]+[\bard])}}{u-v}
    \left(\rule{0ex}{2.1ex}T^{\barb\barc}(u) T^{\bard\bara}(v)-
    T^{\barb\barc}(v) T^{\bard\bara}(u)\right)\\
&=& \frac{(-1)^{[a]'[b]'+([a]'+[b]')[c]'}}{u-v} 
\left(\rule{0ex}{2.4ex}\wt{T}^{cb}(u) \wt{T}^{ad}(v)-\wt{T}^{cb}(v) 
\wt{T}^{ad}(u)\right)
\eeano 
which is the correct expression for the commutator in  $Y(N|M)$, 
since $[\ ]'$ is the correct gradation of $Y(N|M)$.

For the Hopf structure, one has:
\beano
\Delta_{M|N}\, \wt{T}^{ab}(u) &=&
\sum_{c=1}^{K}\, (-1)^{[\bara]([\barb]+1)}
{T}^{\barb c}(u)\otimes{T}^{c\bara}(u)\\ 
&=& \sum_{c=1}^{K}\, 
(-1)^{[\bara]([\barb]+1)+[\barc]([\barb]+1)+[\bara]([\barc]+1)}
\wt{T}^{cb}(u)\otimes\wt{T}^{ac}(u)\\
&=&
P\, \wt{T}^{ac}(u)\otimes\wt{T}^{cb}(u)\, P 
= \Delta'_{N|M}\, \wt{T}^{ab}(u)
\eeano
where we have denoted by $\Delta_{M|N}$ (resp. $\Delta_{N|M}$) the 
coproduct on $Y(M|N)$ (resp. on $Y(N|M)$).

A simple calculation shows
\be
\eps_{M|N}(\wt{T}^{ab}(u)) = \eps'_{N|M}(\wt{T}^{ab}(u))
\ee
where $\eps'=\eps$ is the counit associated to $\Delta'$,
while for the antipode
\be
S'_{M|N}(\wt{T}^{ab}(u)) = 
(-1)^{[\bara]([\barb]+1)}\theta_{a}\theta_{b} 
\left(T^{-1}(u)\right)^{\bara\barb} = S'_{N|M}(\wt{T}^{ab}(u)) =
\left(\wt{T}^{-1}{}'(u)\right)^{ab}
\ee
where in the last equality, the inverse $\wt{T}^{-1}{}'(u)$ is computed 
using $m'$ instead of $m$:
\be
m'\left( \wt{T}^{ab}(u)\otimes \wt{T}^{-1}{}'(u)^{bc}\right)=
(-1)^{([a']+[b'])([b']+[c'])}\, \wt{T}^{-1}{}'(u)^{bc}\cdot \wt{T}^{ab}(u)
=\delta^{ac}
\ee
$S'_{M|N}$ obeys the relations
\be
m'(S'\otimes id)\Delta'=\eps'=m'(id\otimes S')\Delta'
\ee
\finprf

\subsection{Finite dimensional irreducible representations of 
$Y(M|2n)$}
The finite dimensional irreducible representations of the superYangian 
$Y(M|N)$ have been studied in \cite{zhang}. We recall here its main 
results, refering to \cite{zhang} for the proofs. 

We will specify to the case $N=2n$, for it is the only case that is 
needed for twisted superYangians. Moreover, in order to be able to 
deal with the twisted superYangians, we need to choose a positive 
roots system different from the one chosen in \cite{zhang}. Indeed, 
the situation in analogous to the one encountered in the case of 
simple Lie superalgebras, which admit different {\it inequivalent} 
systems of simple roots (see for instance \cite{dico}). For our 
purpose, we define:
\begin{defi}[Positive roots]\hfill\\
Let $\Phi^{\pm,0}\in\NN_{K}^{2}$, where 
$\NN_{K}=\NN_{M+2n}=[1,M+2n]\cap\NN$, 
be defined by
\be
\Phi^{+}=\Big\{(a,b)\in \NN_{K}^{2}, \mbox{ with either }
\left|\begin{array}{l}
1\leq a<b\leq M \\
M+1\leq a<b\leq M+2n \\
1\leq a\leq M,\, M+n+1\leq b\leq M+2n\\
M+n+1\leq a\leq M+2n,\, 1\leq b\leq M
\end{array}\right.
\Big\}
\ee
\be
\Phi^{-}=\{(a,b)\in \NN_{M+2n}^{2}, \mbox{ such that }
(b,a)\in\Phi^{+}\}
\ee
\be
\Phi^{0}=\{(a,a),\mbox{ with }a\in \NN_{M+2n}\}
\ee
We have $\NN_{M+2n}=\Phi^-+\Phi^{0}+\Phi^{+}$, and the positive 
roots will be associated to $T^{ab}(u)$ with $(a,b)\in\Phi^{+}$.
\end{defi}
Once the positive roots system is chosen, one can introduce the 
notion of highest weight vectors:
\begin{defi}[Highest weight vectors]\hfill\\
Let $\cm$ be a module of $Y(M|2n)$. A highest weight vector $v\in\cm$ 
is defined by
\be
\begin{array}{l}
T^{aa}(u)v=\lbda^a(u)v,\ \forall a=1,\ldots,M+2n\\
T^{aa}(u)v=0,\ \forall (a,b)\in\Phi^{+}
\end{array}
\ee
$\lbda(u)=(\lbda^{1}(u),\ldots,\lbda^{M+2n})\in\CC[[u^{-1}]]$ is the 
highest weight associated to $v$.
\end{defi}
The notion of highest weight vectors grounds in the to following 
properties:
\begin{prop}
Any irreducible finite dimensional representation of $Y(M|N)$ admits 
a unique (up to multiplication by scalars) highest weight vector.
\end{prop}
\begin{prop}\label{Irreps-YMN}
The irreducible representation of $Y(M|N)$ with highest weight 
$\lbda(u)$ is finite dimensional if and only if we have:
\bea
\frac{\lbda^a(u)}{\lbda^{a+1}(u)} &=& \frac{P_{a}(u+1)}{P_{a}(u)},\ 
1\leq a\leq M+N-1,\ a\neq M \label{eq1-irrepsYMN}\\
\frac{\lbda^M(u)}{\lbda^{M+1}(u)} &=& \frac{P_{M}(u)}{{P}_{M+N}(u)}
\label{eq2-irrepsYMN}
\eea
where $P_{a}(u)$ are monic polynomials.
\end{prop}

Let us remark that, some signs differ between our 
presentation and the presentation given in \cite{zhang} because of the 
definition for $T(u)$: the relation between these two notations is 
given by $T^{ab}_{(n)}=(-1)^{[b]}t_b^a[n]$.
\begin{defi}[Evaluation representations]\hfill\\
Let $J^{ab}$ be the generators of the $gl(M|N)$ superalgebra and
$\pi^{ab}=\pi(J^{ab})$ a finite dimensional representation of 
$gl(M|N)$. Then, the morphism
\be
ev(T(u))=1+\frac{\EE}{u}\mb{with} \EE=\pi^{ab}E_{ab}\label{eval}
\ee 
provides a finite dimensional representation of $Y(M|N)$, called an evaluation 
representation.
\end{defi}
The usefulness of evaluation representations reveals in the following 
theorem:
\begin{theor}
Any irreducible finite dimensional representation of $Y(M|N)$ is 
isomorphic to the irreducible part of tensor products of evaluation 
representations.
\end{theor}

\sect{Twisted superYangians\label{superYtw}}
\subsection{Introduction to $Y(M|N)^\tau$}
We now introduce the notion of twisted super Yangian, in the same way 
twisted Yangians have been defined from the Yangians $Y(N)$.

We first introduce the transposition $t$ on matrices:
\be
\begin{array}{l}
E^t_{ab}=(-1)^{[a]([b]+1)}\theta_{a}\theta_{b}\, E_{\barb\bara}
\ \mbox{ with }\ \theta_{a}=\pm1
\\[1.2ex]
\bara=M+1-a\ \mbox{ for }1\leq a\leq M\\[1.2ex]
\bara=2M+N+1-a\ \mbox{ for }M+1\leq a\leq M+N
\end{array}
\label{transp}
\ee
which satisfies 
\be
(AB)^t=B^t\, A^t
\ee
Demanding the transposition to be of order 2 leads to the constraint
\be
(-1)^{[a]}\, \theta_{a}\theta_{\bara}=\theta_{0}=\pm1\ \forall\ a
\ee
Let us stress that $\bara$ has \underline{not} the same meaning as 
in section \ref{supY}: from 
now on, we will use this notation to denote (\ref{transp}). These 
new $\bara$ satisfy $[a]=[\bara]$.

Let us also note that there is a freedom on the definition of the 
transposition, 
$\theta_{a}\rightarrow (-1)^{[a]}\theta_{a}$: this freedom is fixed 
when imposing $\II^t=\II$.
Remark also the useful identity
\be
(-1)^{[a]([b]+1)}\theta_{a}\theta_{b}=
(-1)^{[\bara]([\barb]+1)}\theta_{\bara}\theta_{\barb}
\ee

Then, we define on $Y(M|N)$:
\be
\tau[T(u)]=\sum_{{a,b}}\tau[T^{ab}(u)]\, E_{ab}=\sum_{{a,b}}T^{ab}(-u)\, 
E^t_{ab}
\ee
which reads for the superYangian generators:
\be
\tau(T^{ab}(u))=(-1)^{[a]([b]+1)}\theta_{a}\theta_{b}\ T^{\barb\bara}(-u)
\label{supertau}
\ee

\begin{prop}
    $\tau$ is  an algebra automorphism for $Y(M|2n)$ and 
    $Y(2m|N)$ only. 
    
    In that case, one must choose $\theta_{0}=+1$ for $Y(M|2n)$ and 
    $\theta_{0}=-1$ for $Y(2m|N)$.
\end{prop}
\prf
Considering   subalgebras mentioned in property \ref{subalg}, 
one can see that $\tau$ acts as an automorphism of $Y(M)$ of the type 
defined in \cite{rag}, with $\theta_{a}\theta_{\bara}=\theta_{0}$, and 
an automorphism of $Y(N)$ with $\theta_{a}\theta_{\bara}=-\theta_{0}$.
Using the results of 
\cite{rag}, where it is proved that when $M$ is odd, 
one must have $\theta_{0}=+1$ in $Y(M)$, one immediately concludes 
that we cannot have $MN$ odd, and that the values for $\theta_{0}$ are 
the ones given in the property. 

Then, it is a simple matter of calculation to show that $\tau$ is an 
automorphism of the superalgebra $Y(M|N)$:
\be
\tau\left([T^{ab}(u),T^{cd}(v)\}\right) = 
[\tau(T^{ab}(u)),\tau(T^{cd}(v))\}
\ee
\finprf
One defines in $Y(M|N)$ (we take $MN$ even):
\bea
S(u) &=& T(u)\, \tau[T(u)]
=\sum_{a,b=1}^{M+N}S^{ab}(u)E_{ab}
=\II+\sum_{a,b=1}^{M+N}\sum_{n>0}u^{-n}S^{ab}_{(n)}\\
S^{ab}_{(n)} &=& \sum_{c=1}^{M+N}\sum_{p=0}^n(-1)^{p}
(-1)^{[c]([b]+1)} 
\theta_{c}\theta_{b}T^{ac}_{(n-p)}T^{\barb\barc}_{(p)}
\label{Sabn}\\
S^{ab}(u) &=& \sum_{c=1}^{M+N}
(-1)^{[c]([b]+1)} 
\theta_{c}\theta_{b}T^{ac}(u)T^{\barb\barc}(-u)
\eea
$S(u)$ defines a subalgebra of the superYangian:
\begin{theor}\label{grY}
    $S(u)$ obey the following relations:
\be
R_{12}(u-v)\, S_{1}(u)\, R'_{12}(u+v)\, S_{2}(v) = 
S_{2}(v)\, R'_{12}(u+v)\, S_{1}(u)\, R_{12}(u-v)
\label{rsrs}
\ee
\be
\tau(S(u))=S(u)+\frac{\theta_{0}}{2u}(S(u)-S(-u))
\label{tauS}
\ee
where $R(x)$ is the superYangian $R$-matrix,
\be
R'(x)=\II+\frac{1}{x}Q=R^{t_{1}}(-x) \ \mbox{ with }\ Q=P^{t_{1}}\label{matRt}
\ee
 and $t_{1}$ is the transposition 
(\ref{transp}) in the first space.
These two relations \underline{uniquely} define a subalgebra $Y(M|N)^\tau$ in the 
super Yangian.
\end{theor}
\prf
One starts with the relation (\ref{RTT}), applies the transposition 
$t_{1}$ and the sign operation $(u,v)\rightarrow (-u,-v)$ to get
\be
\tau[T_{1}(u)]R'_{12}(u+v)T_{2}(v) = T_{2}(v)R'_{12}(u+v)\tau[T_{1}(u)]
\label{tauTT}
\ee
A direct calculation shows that
\be
P\, Q =Q\,P=\theta_{0}Q\ ;\ P^2=\II\ \mbox{ and }\ Q^{2}=(M-N)Q
\label{propQ}
\ee
Thus, applying $P(.)P$ on (\ref{tauTT}) leads to (after the exchange 
$u\leftrightarrow v$):
\be
T_{1}(u)R'_{12}(u+v)\tau[T_{2}(v)] = \tau[T_{2}(v)]R'_{12}(u+v)T_{1}(u)
\label{TtauT}
\ee
Finally, applying once again the transposition $t_{1}$ and 
$(u,v)\leftrightarrow(-u,-v)$, we obtain:
\be
R_{12}(u-v)\tau[T_{1}(u)]\tau[T_{2}(v)] = 
\tau[T_{2}(v)]\tau[T_{1}(u)]R_{12}(u-v)
\label{tauTtauT}
\ee
which is another way to prove that $\tau$ is an automorphism. A simple 
calculation using (\ref{RTT}), (\ref{tauTT}), (\ref{TtauT}) and 
(\ref{tauTtauT}) shows then that (\ref{rsrs}) is satisfied. 

The second relation is also proved directly:
\bea
\left(\tau[S(u)]\right)^{ab} &=& \sum_{c=1}^{M+N} 
(-1)^{[a]([b]+1)+[c]([a]+1)} \theta_{a}\theta_{b}\theta_{\bara}\theta_{c}
T^{\barb c}(-u)T^{a\barc}(u)\\
&=& S(u)+\sum_{c=1}^{M+N} 
(-1)^{[a]([b]+[c])} \theta_{b}\theta_{c}
[T^{\barb c}(-u),T^{a\barc}(u)\}\\
&=& \left(S(u)+\frac{\theta_{0}}{2u}(S(u)-S(-u))\right)^{ab}
\eea
where, in the last step, we have used the graded commutator (\ref{comYMN}).

\null

Conversely, let us start with a subalgebra $\ca$ of $Y(M|2n)$ whose 
generators $\sigma_{(n)}^{ab}$ obey (\ref{rsrs}) and (\ref{tauS}).
There is an obvious surjective morphism $\jmath$ from $\ca$ to
$Y(M|2n)^{+}$. Thus, it remains to show that this morphism is 
injective. We follow the argumentation done in \cite{Molev} for the case 
of (non super) twisted Yangians.

We first introduce a filtration on $Y(M|2n)$ induced by
\be
\dg(T_{(p)}^{ab})=p \mb{and} \dg(XY)=\dg(X)\dg(Y),\ \forall X,Y\in Y(M|2n)
\ee
The graded algebra $gr Y(M|2n)$ is defined as usual by:
\bea
Y_{p}(M|2n) &=& \{X\in Y(M|2n), \mb{with} \dg(X)\leq p\}\ ;\ p>0\\
Y_{0}(M|2n) &=& \CC\ ;\ gr_{0}Y(M|2n)=\CC\\
gr_{p}Y(M|2n) &=& Y_{p}(M|2n)\, /\, Y_{p-1}(M|2n)\ ;\ p>0\\
grY(M|2n) &=& \oplus_{p\geq0}gr_{p}Y(M|2n)
\eea
Since for $X\in gr_{p}Y(M|2n)$ and $Y\in gr_{q}Y(M|2n)$, we have 
$[X,Y\}\in gr_{p+q-1}Y(M|2n)$, we deduce that $gr Y(M|2n)$ is 
commutative\footnote{During this proof, and to avoid confusion with 
the gradation deg, we will write commutative and commutator where one 
should has written $\ZZ_{2}$-graded commutative and 
$\ZZ_{2}$-graded commutator.}. The same is true for $gr Y(M|2n)^{+}$, 
here the filtration is induced by the $Y(M|2n)$ one.

Similarly, on $\ca$, we introduce a filtration given by
\be
\dg(\sigma_{(n)}^{ab})=n
\ee
This makes $gr\ca$ a commutative algebra, for the same reasons as 
above. Moreover, since the morphism $\jmath$ preserves the 
filtration, it is enough to show that the induced morphism 
$\bar\jmath$ between graded algebras is injective.

Let $\bar T_{(p)}^{ab}$ and $\bar S_{(p)}^{ab}$ be the image of 
$T_{(p)}^{ab}$ and $S_{(p)}^{ab}$ in $gr Y(M|2n)^{+}$. The 
expression (\ref{Sabn}) is still valid for the elements of $gr 
Y(M|2n)^{+}$, so that we deduce
\be
\bar S_{(p)}^{ab}= (-1)^p (-1)^{[a]([b]+1)}\theta_{a} \theta_{b}  
\bar S_{(p)}^{\barb\bara}\label{tauSbar}
\ee
Thus, we conclude that $gr Y(M|2n)^{+}$ is isomorphic to the algebra 
of polynomials in the ($\ZZ_{2}$-graded) letters $x_{(p)}^{ab}$ 
submitted to the constraints (\ref{tauSbar}). 

Finally, the symmetry relation (\ref{tauS}) in the algebra $gr\ca$ just
takes the form (\ref{tauSbar}), so that $gr\ca$ is also isomorphic to 
the algebra of polynomials in the ($\ZZ_{2}$-graded) letters $x_{(p)}^{ab}$ 
submitted to the constraints (\ref{tauSbar}). Hence, $\bar\jmath$ is 
injective.
\finprf

\begin{coro}[PBW basis for $Y(M|2n)^{\tau}$]\hfill\\
\label{PBW}
Given an arbitrary linear order on the following set of generators 
(for $n=1,2,\ldots$):
\be
\begin{array}{lll}
S^{ij}_{(2n)} &\mbox{for } 1\leq i,j \leq  M
&\mbox{and }\ i+j\leq M+1\\
S^{ij}_{(2n+1)} &\mbox{for } 1\leq i,j \leq M
&\mbox{and }\ i+j< M+1\\
S^{ij}_{(2n)} &\mbox{for } M+1\leq i,j \leq  M+2n
&\mbox{and }\ i+j< 2M+2+2n\\
S^{ij}_{(2n+1)} &\mbox{for } M+1\leq i,j \leq M+2n
&\mbox{and }\ i+j\leq 2M+2+2n\\
S^{ij}_{(2n)} &\mbox{for } M+1\leq i \leq  M+2n
&\mbox{and }\ 1\leq j\leq M\\
\end{array}
\ee
any element of $Y(M|2n)^{\tau}$ is uniquely written as a linear 
combination of the ordered monomials in these generators.
\end{coro}
\prf
It is a direct consequence of the proof of theorem \ref{grY}. Indeed, 
considering $grY(M|2n)^{\tau}$, it is sufficient to find a basis for 
it, i.e. for the algebra of polynomials in the ($\ZZ_{2}$-graded) 
letters $x_{(p)}^{ab}$ submitted to the constraints (\ref{tauSbar}).
From the property $a+b\leq M+1\ \Leftrightarrow\ \bara+\barb\geq M+1$, 
$\bara=b\ \Leftrightarrow\ a+b= M+1$ when $a,b\leq M$, and 
$a+b\leq M+n+2\ \Leftrightarrow\ \bara+\barb\geq M+n+2$, 
$\bara=b\ \Leftrightarrow\ a+b= M+n+2$ when $a,b\geq M+1$,
an analysis of these constraints lead to the above basis.
\finprf

Although several automorphisms $\tau$ can be defined (depending upon the 
choices for the $\theta$'s), they all lead to the same subalgebra $Y(M|N)^\tau$:
\begin{prop}
All the $\theta_{a}$ dependence can be removed in the commutation 
relations of  $Y(M|2n)^\tau$.
\end{prop}
\prf
We prove the  property by exhibiting a basis in which the $\theta$ 
dependence has disappeared. When restricted to the bosonic part, it is 
the same basis as the one given in \cite{rag} for (bosonic) twisted 
Yangians.

For $Y^-(2n)$:
\be
{J}^{ij}(u)=\theta^{i}\theta^{j}S^{ij}(u)\ ;\  
{K}^{ij}(u)= \theta^{i}S^{i,\barj}(u)\ ;\  
{\bar K}^{ij}(u)= \theta^{j}S^{\bari,j}(u)\ ;\ i,j=1,\ldots,n
\label{bosN}
\ee
For $Y^+(M)$, the redefinition is the same as before, plus for 
the remaining generators which appear when $M=2m+1$:
\be
{J}_{0}(u)=\theta^{\tilm}S^{\tilm\tilm}(u)\ ;\ {L}^{i}(u)=\theta^{i} 
S^{\tilm,i}(u)\ ;\ 
{\bar L}^{i}(u)=\theta^{\tilm}\theta^{i} S^{i,\tilm}(u)\ ;\ i,j=1,\ldots,m
\label{bosM}
\ee
This proves that the commutation relations among generators of 
$Y^-(2n)$, and those among $Y^+(M)$ are free from $\theta$'s in this 
basis. Commuting an element of $Y^-(2n)$ with one of $Y^+(M)$ provide 
the change of basis for the fermionic generators:
\be
\begin{array}{lll}
F^{ai}(u)=\theta_{a}\theta_{i}S^{ai}(u) & 
\bar F^{ia}(u)=\theta_{a}\theta_{i}S^{ai}(u) & i=M+1,\ldots,M+n\\
G^{ai}(u)=\theta_{a}S^{a\bari}(u) & 
\bar G^{ia}(u)=\theta_{i}S^{\bara i}(u) & a=1,\ldots,M\\
H^{i}(u)=\theta_{i}S^{m+1,i}(u) & 
\bar H^{i}(u)=\theta_{m+1}\theta_{i}S^{i,m+1}(u) & \mbox{if } M=2m+1
\end{array}
\label{ferm}
\ee
All the $Y(M|2n)^+$-generators are expressible in terms of the generators 
(\ref{bosN}), (\ref{bosM}) and (\ref{ferm}) using the symmetry 
relation (\ref{tauS}).
Then, one can check that all the graded commutators in this basis are 
free from $\theta$.
\finprf
\begin{defi}
    The twisted superYangian $Y(M|2n)^+\equiv Y(2n|M)^-$ is the 
    subalgebra generated 
    by $S(u)=T(u)\tau[T(u)]$, with $\tau$ given in (\ref{supertau}) 
    and 
    \be\begin{array}{l}
    \theta_{a}=1\ \mbox{ for }1\leq a\leq M\\
    \theta_{a}=\mbox{sg}(\frac{2M+2n+1}{2}-a)\ \mbox{ for }M+1\leq a\leq M+2n
    \end{array}
 \ee
 \end{defi}
 \subsection{Few properties of $Y(M|2n)^+$}
  \begin{prop}
The relation (\ref{rsrs}) is equivalent to the 
following commutator:
 \bea
{[}S_{1}(u),S_{2}(v)] &=& \frac{1}{u-v}\left( \rule{0ex}{2.4ex}
P_{12}S_{1}(u)S_{2}(v) - S_{2}(v)S_{1}(u)P_{12}\right)+ \nonu
&&-\frac{1}{u+v}\left( \rule{0ex}{2.4ex}
S_{1}(u)Q_{12}S_{2}(v)-S_{2}(v)Q_{12}S_{1}(u)
\right)+ \label{comS1S2}\\
 &&+\frac{1}{u^2-v^2}\left(\rule{0ex}{2.4ex}
P_{12}S_{1}(u)Q_{12}S_{2}(v)-S_{2}(v)Q_{12}S_{1}(u)
P_{12}\right)\nonumber
\eea
and also to
 \bea
&&{[}S^{ab}(u),S^{cd}(v)\} = \frac{(-1)^{([a]+[b])[c]}}{u-v}\,
(-1)^{[a][b]}\left( \rule{0ex}{2.4ex}
S^{cb}(u)S^{ad}(v) - S^{cb}(v)S^{ad}(u)\right)+ \nonu
&&\hspace{2.1ex}-\frac{(-1)^{([a]+[b])[c]}}{u+v}\left( \rule{0ex}{2.4ex}
(-1)^{[a][c]}\theta_{b}\theta_{\barc}S^{a\barc}(u)S^{\barb d}(v)
-(-1)^{[b][d]}\theta_{\bara}\theta_{d}S^{c\bara}(v)S^{\bard b}(u)
\right)+ \nonu
&&\hspace{2.1ex} +\frac{(-1)^{([a]+[b])[c]}}{u^2-v^2}\,
(-1)^{[a]}\theta_{a}\theta_{b}\left(\rule{0ex}{2.4ex}
S^{c\bara}(u)S^{\barb d}(v)-S^{c\bara}(v)S^{\barb d}(u) \right)
\label{comSijSkl}
\eea
 \end{prop}
 \prf
 (\ref{comS1S2}) follows from a direct calculation using 
 (\ref{rsrs}), (\ref{matR}) and (\ref{matRt}). 
 \finprf
As an obvious consequence, we get
 \begin{coro}\label{ospMN}
 The twisted super-Yangian $Y(M|2n)^+$ contains $osp(M|2n)$ as Lie 
 sub-superalgebra. It is generated by
\be
S^{ab}_{(1)}= T^{ab}_{(1)}-(-1)^{[a]([b]+1)}\theta_{a}\theta_{b}\ 
T^{\barb\bara}_{(1)}
\;\ a,b=1,\ldots,M+2n
\ee
which obey
\be
S^{ab}_{(1)}=-(-1)^{[a]([b]+1)}\theta_{a}\theta_{b}\ S^{\barb\bara}_{(1)}
\ee
and defines a morphism of algebra $\cu[osp(M|2n)]\rightarrow Y(M|2n)^+$.

The action of the $osp(M|2n)$ generators on the twisted Yangian is given by:
 \bea
{[}S^{ab}_{(1)},S^{cd}(v)\} &=& (-1)^{([a]+[b])[c]}\left\{\,
(-1)^{[a][b]}\left( \rule{0ex}{2.4ex}
\delta_{cb}S^{ad}(v) - \delta_{ad}S^{cb}(v)\right)\right.+ \nonu
&&\hspace{2.1ex}-\theta_{\bara}\theta_{b}\left.\left( \rule{0ex}{2.4ex}
\delta_{a\barc}S^{\barb d}(v)-\delta_{\bard b}S^{c\bara}(v)
\right) \right\}
\label{comJSabcd}
\eea
\end{coro}
 \prf
 Expanding $(u\pm v)^{-1}=u^{-1}(1\mp v u^{-1}+\ldots)$ and taking the coefficient 
 of $u^{-1}v^{-1}$ in (\ref{comS1S2}) leads to:
\be
{[}S_{1(1)},S_{2(1)}\} = P_{12}S_{2(1)} - S_{2(1)}P_{12}-
Q_{12}S_{2(1)}+S_{2(1)}Q_{12} \label{comS1}
\ee
where the subscript $(1)$ refers to the coefficient of $u^{-1}$ and 
$v^{-1}$, while the indices $1,2$ label the auxiliary spaces. The symmetry 
relation projected on the $u^{-1}$ term reads
\be
S_{(1)}^t=-S_{(1)} \label{symS1}
\ee
(\ref{comS1}) and (\ref{symS1}) 
are just the defining relations of $osp(M|2n)$.

Starting now from (\ref{comSijSkl}) and taking the coefficient of $u^{-1}$ 
gives
the relation (\ref{comJSabcd}). Note that taking the coefficient of $v^{-1}$
in this last relation gives again the commutation relations of $osp(M|2n)$.
\finprf

 Let us denote the $osp(M|2n)$ generators by $J^{ab}$ and gather 
 them in the matrix
 \be
F=\sum_{a,b=1}^{M+N}J^{ab} F_{ab}
\mb{with}
F_{ab}=E_{ab}-(-1)^{[a]([b]+1)}\theta_{a}\theta_{b}\ E_{\barb\bara}
 \ee
 It satisfies:
 \be
 \begin{array}{l}
 F^t=-F \\[1.2ex]
 {[}F_{1},F_{2}\}=P_{12}F_{2}-F_{2}P_{12}+F_{2}Q_{12}-Q_{12}F_{2}
\end{array}\label{comOsp}
 \ee
 where $t$ is the transposition (\ref{transp}).
\begin{prop}
 The following map defines an algebra inclusion:
 \be
 \begin{array}{l}
Y(M|2n)^+\ \rightarrow\ \cu[osp(M|2n)]\\[1.2ex]
\displaystyle
S(u) \ \rightarrow\ \FF(u)=\II+\frac{1}{u+\half}F
\end{array}
\ee
\end{prop}
\prf
We have to prove that $\FF$ obeys to the relations (\ref{rsrs}) and 
(\ref{tauS}). A direct calculation shows:
\[
\FF^t(-u)=\II+\frac{1}{u-\half}F=\FF(u)+\frac{1}{2u}\left(\FF(u)-\FF(-u)\right)
\]
Moreover, using (\ref{propQ}), the commutator (\ref{comOsp}) and the 
relations 
 \be
\begin{array}{l}
P_{12}F_{2}=F_{1}P_{12}\ \Rightarrow\ Q_{12}F_{2}=-Q_{12}F_{1}\\
P_{12}F_{1}=F_{2}P_{12}\ \Rightarrow\ F_{1}Q_{12}=-F_{2}Q_{12}
\end{array}\label{QS1}
\ee
one proves that we have 
\beano
{[}\FF_{1}(u),\FF_{2}(v)\} &=& \frac{1}{u-v}\left( \rule{0ex}{2.4ex}
P_{12}\FF_{1}(u)\FF_{2}(v)
-\FF_{2}(v)\FF_{1}(u)P_{12}\right)+\\
&&-\frac{1}{u+v}\left( \rule{0ex}{2.4ex}
\FF_{1}(u)Q_{12}\FF_{2}(v)
-\FF_{2}(v)Q_{12}\FF_{1}(u)\right) +\\
&&+ \frac{1}{u^2-v^2}\left( \rule{0ex}{2.4ex}
P_{12}\FF_{1}(u)Q_{12}\FF_{2}(v)
-\FF_{2}(v)Q_{12}\FF_{1}(u)P_{12}\right) 
\eeano
\finprf
The relation between $osp(M|2n)$ and $Y(M|2n)^+$ also reveals in
\begin{prop}
    $Y(M|2n)^+$ is a deformation of ${\cu}\left(osp(M|2n)[x]\right)$, 
    the (positive modes) loop algebra based on $osp(M|2n)$.
 \end{prop}
 \prf
 We start with $S(u)=\II+s(u)$, and make a change of basis 
 $\wt{s}(u)=\hbar^{-1}\, s(u/\hbar )$. In this basis, the commutation 
relations 
 read:
\beano
{[}\wt{s}_{1}(u),\wt{s}_{2}(v)\} &=& \frac{1}{u-v}\left( \rule{0ex}{2.4ex}
 P_{12}\wt{s}_{1}(u)+P_{12}\wt{s}_{2}(v) - \wt{s}_{1}(u)P_{12}- 
\wt{s}_{2}(v)P_{12}\right)
 + \\
&&-\frac{1}{u+v}\left( \rule{0ex}{2.4ex}
 \wt{s}_{1}(u)Q_{12}+Q_{12}\wt{s}_{2}(v)-Q_{12}\wt{s}_{1}(u)-
\wt{s}_{2}(v)Q_{12}
\right)+\\
 && +\frac{\hbar}{u-v}\left( \rule{0ex}{2.4ex}
P_{12}\wt{s}_{1}(u)\wt{s}_{2}(v) - 
 \wt{s}_{2}(v)\wt{s}_{1}(u)P_{12}\right)+\\
 &&-\frac{\hbar}{u+v}\left( \rule{0ex}{2.4ex}
 \wt{s}_{1}(u)Q_{12}\wt{s}_{2}(v)-\wt{s}_{2}(v)Q_{12}\wt{s}_{1}(u)
 \right)+ \\
 &&+\frac{\hbar}{u^2-v^2}\left(\rule{0ex}{2.4ex}
 P_{12}\wt{s}_{1}(u)Q_{12}+P_{12}\wt{s}_{2}(v)-Q_{12}\wt{s}_{1}(u)
 P_{12}-\wt{s}_{2}(v)P_{12}\right)+\\
 &&+\frac{\hbar^2}{u^2-v^2}\left(\rule{0ex}{2.4ex}
 P_{12}\wt{s}_{1}(u)Q_{12}\wt{s}_{2}(v)-\wt{s}_{2}(v)Q_{12}\wt{s}_{1}(u)
 P_{12}\right)\\
\eeano
For $\hbar\neq0$ all the algebras $Y_{\hbar}(M|2n)^+$ are isomorphic, 
while in the limit $\hbar\rightarrow 0$, $Y_{\hbar=0}(M|2n)^+$ 
reduces to 
 $\cu\left(osp(M|2n)[x]\right)$.
 \finprf
Note also the isomorphism
\begin{prop}[Automorphism of $Y((M|2n)^+$]\hfill\\
The transformations
\be
S(u)\rightarrow g(u)\, S(u)\mb{with $g(u)$ even 
$\CC$-function}\label{g(u)}\\
\ee
are automorphisms of $Y(M|2n)^+$.
\end{prop}
\prf
Multiplying (\ref{rsrs}) by $g(u)g(v)$ shows that it is invariant 
under the transformation (\ref{g(u)}), for any function $g$. 
The symmetry relation (\ref{tauS}) is preserved for $g(u)$ even only.
\finprf
There is another of type of automorphism that we will need when looking at the representations of twisted superYangians:
\begin{defi}[$\#$ involution]\hfill\\
\label{defhash}
For any index $i=1,...,M+N$, we define $i'$ by
\be
i'=\left\{\begin{array}{l} m'=\tilm+1,\ \ (\tilm+1)'=m\\
i\mb{otherwise}
\end{array}\right.
\ee
 $\#$ defined by
\be
S_\#^{ij}(u)=S^{i'j'}(u)
\ee
is an order 2 automorphism of $Y(M|2n)^+$.
\end{defi}
\prf
Obvious direct calculation from the relations (\ref{comSijSkl}) and 
(\ref{tauS}).
\finprf

 For the Hopf structure, and mimicking again the case of twisted Yangians, 
one can show:
 \begin{prop}
     $Y(M|2n)^+$ is a left coideal of $Y(M|2n)$: 
\be
\Delta(Y(M|2n)^+)\subset Y(M|2n)\otimes Y(M|2n)^+
\ee
More precisely:
 \beano
\Delta(S^{ab}_{(p)}) &=& 
\sum_{y=0}^{p}\sum_{q=0}^{y}\sum_{d,e=1}^{M+2n}\ 
(-1)^q(-1)^{[d]([e]+[b])}\theta_{b}\theta_{e}\ 
T^{ad}_{(y-q)}T^{\barb\bare}_{(q)}\otimes S^{de}_{(y)}\\
\Delta[S^{ab}(u)] &=& \sum_{d,e=1}^{M+2n}\ 
(-1)^{[d]([e]+[b])}\theta_{b}\theta_{e}\ 
 T^{ad}(u)T^{\barb\bare}(-u)\otimes S^{de}(u)
\eeano
 \end{prop}
 \prf
 Direct calculation using (\ref{DeltaT}) and (\ref{Sabn}).
 \finprf
\subsection{Subalgebras of $Y(M|2n)^+$}
 \begin{prop}
 The twisted superYangian $Y(M|N)^+$ contains as subalgebras 
 $Y(M)^{+}$, $Y(N)^-$ and $osp(M|N)$.
 \end{prop}
 \prf
 A direct examination on the commutator (\ref{comYMN}) and the symmetry 
 relation (\ref{tauS}) shows that $S^{ab}(u)$ with $a,b=1,\ldots,M$ 
 (resp. $a,b=M+1,\ldots,M+N$) generates the twisted Yangian $Y(M)^{+}$
 (resp. $Y(N)^-$). The last inclusion has been proved in the 
 corollary \ref{ospMN}.
 \finprf
 \begin{prop}
 As algebra embbedings, we have:
 \be
 Y(1|2)^+\subset Y(2m+1|2n)^+ \mb{and} Y(2|2)^+\subset Y(2m|2n)^+
 \ee
 \end{prop}
 \prf
 We consider the generators $S^{ij}(u)$, with $i,j=m+1,2m+n+1, 
 2m+n+2$ in $Y(2m+1|2n)^+$ and $S^{ij}(u)$, with $i,j=m,m+1,2m+n, 
 2m+n+1$ in $Y(2m|2n)^-$: they obey the commutation and symmetry 
 relations of $Y(1|2)^+$ and $Y(2|2)^+$ respectively.
 \finprf
 
 Note that there is no 
 regular embbeding of $Y(1|2)^+$ into $Y(2m|2n)^+$. The circumstances 
 are 
 here different from both simple superalgebras and non-super twisted 
 Yangians cases: in the first case, it 
 always exists a regular $osp(1|2)$ embbeding, and in the second case, 
 one can always construct a regular $Y(2)^\pm$ embbeding in the 
 twisted Yangian $Y(M)^\pm$. It is the symmetry relation which causes 
 this unusual situation. 
 \begin{prop}
 As algebra embbeding, we have:
 \be
 Y(2m|2n)^+\subset Y(2m+1|2n)^+ 
 \ee
Let us stress however that $ Y(2m|2n)^+$ is \underline{not} a
Hopf coideal of $ Y(2m+1|2n)^+$.

The same results apply for $Y(2m)^+$ and $ Y(2m+1)^+$.
 \end{prop}
 \prf
 Let $s^{ij}(u)$ be the generators of $Y(2m+1|2n)^+$. 
 We set $M=2m+1$ and introduce 
 \[
 \begin{array}{ll}
 \sigma^{ij}(u)=s^{ij}(u) &\mbox{for }1\leq i,j\leq m
 \mbox{ and }M+1\leq i,j\leq M+n+1
\\[1.2ex]
\sigma_{i,j}(u)=s^{i-1,j-1}(u) &\mbox{for }m+2\leq i,j\leq M
\mbox{ and }M+n+2\leq i,j\leq M+2n
\\[2.1ex]
\sigma^{ij}(u)=s^{i-1,j}(u) &\mbox{for }\left\{
\begin{array}{l}1\leq j\leq m \mbox{ or }M+1\leq j\leq M+n+1 \\ 
m+2\leq i\leq M\mbox{ or }M+n+2\leq i\leq M+2n\end{array}
\right.\\[2.1ex]
\sigma^{ij}(u)=s^{i,j-1}(u) &\mbox{for }\left\{
\begin{array}{l}1\leq i\leq m \mbox{ or }M+1\leq i\leq M+n+1 \\ 
m+2\leq j\leq M\mbox{ or }M+n+2\leq j\leq M+2n\end{array}
\right.
 \end{array}
 \]
 We prove that the generators $\sigma^{ij}(u)$ generates 
 $Y(2m|2n)^+$. We denote by $x\rightarrow \bar x$ the "bar" operator 
 introduced in
 (\ref{supertau}) for $Y(2m+1|2n)^+$, and by $x\rightarrow \wt x$ this 
 "bar" operator for $Y(2m|2n)^+$. In the same way, we call $\tau$ 
 and $\theta$ (resp. $\wt{\tau}$ and $\wt\theta$) the corresponding 
 operations in $Y(2m+1|2n)^+$ (resp. $Y(2m|2n)^+$).
 It is easy to see that 
 \[
 \begin{array}{lll}
 \bar{\imath}=\wt \imath & \theta_{i}=\wt\theta_{i}&\mbox{for }i\leq 
 m\mbox{ and }M+1\leq i\leq M+n+1\\[1.2ex]
 \overline{\imath-1}=\wt \imath\quad & 
 \theta_{i-1}=\wt\theta_{i}&\mbox{for }
 m+2\leq i\leq M\mbox{ and }M+n+2\leq i,j\leq M+2n
 \end{array}
 \]
 so that the action of $\tau$ on $s(u)$ is equivalent to the action of
 $\wt\tau$ on $\sigma(u)$. It also proves that the symmetry relation 
 of $s(u)$ (coming from $Y(2m+1|2n)^+$) implies the symmetry relation 
 for $\sigma(u)$ (as $Y(2m|2n)^+$ generator).
 
 In the same way, one shows, starting with the commutation relations 
 of $s(u)$, that the commutation relations of $\sigma(u)$ are those of 
 $Y(2m|2n)^+$.

Finally, computing $\Delta\sigma^{ij}(u)$ as it is induced from the 
$Y(2m+1|2n)^+$ coproduct does not lead to the coproduct formula for 
$Y(2m|2n)^+$.
 \finprf

\sect{Finite dimensional irreducible representations of $Y(M|2n)^+$
\label{irreps}}
We study here the finite dimensional irreducible representations of 
$Y(M|2n)^{+}$ starting from $Y(M|2n)$ in the same way those of 
$Y(M)^{\pm}$ have been studied starting from $Y(M)$ \cite{Molev}.

As a short hand writing, we note irreps for irreducible representations.

 \subsection{Generalities}
 \begin{defi}[Highest weight vector]\hfill\\
Let $\cm$ be a module of $Y(M|2n)^{+}$. A nonzero vector $v\in \cm$ 
is called highest weight if it satisfies
\bea
S^{ij}(u)v=0&\mbox{ for } & (i,j)\in\Phi_{+}\label{Sv=0}\\
S^{ii}(u)v=\mu_i(u)v &\mbox{ for }& i=1,\ldots,M+2n
\eea
for some formal series $\mu_i(u)\in 1+u^{-1}\CC[[u^{-1}]]$. 
The set $\mu(u)\equiv(\mu_1(u),\ldots,\mu_{M+2n}(u))$ is 
the highest weight of $\cm$. 
\end{defi}

{\bf Remark 1:} Due to the symmetry relation (\ref{tauS}), some of the 
relations (\ref{Sv=0}) are redundant, and one could reduce $\Phi^{+}$: 
we keep it as it is to make the comparison with the $Y(M|N)$ case. 

Note also that, in the basis of \cite{zhang}, the symmetry relation 
would have led to $S^{ij}(u)v=0$, $\forall i\neq j$, hence the 
present choice for the positive roots system.

\null

{\bf Remark 2:} The symmetry relation also implies for the highest 
weight:
\be
\label{symmetry-mu}
\mu_{\bar{a}}(u)=\frac{1}{2u}\mu_a(u)+\frac{2u-1}{2u}\mu_a(-u)
\ee
so that, in the $Y(2m+1|2n)^+$ case, $\mu_{m+1}(u)$ is an even 
function of $u$.
 \begin{defi}[Highest weight representations]\hfill\\
\label{HWR}
A  representation $V$ of the twisted super Yangian $Y(M|2n)^+$ is called  
highest weight  it is generated by a highest weight vector . 
If $\mu(u)$ is the highest weight of $v$, we will use the notation
$V[\mu(u)]$ for $V$.
\end{defi}

\begin{theor}
Every finite-dimensional irrep $V$ of $Y(M|2n)^+$ is 
highest weight. Moreover, V contains a unique (up to scalar multiples) 
highest weight vector.
\end{theor}

\prf 
We define 
\be
V_+=\{v\in V|S^{ab}_{(p)}v=0,\;\forall(a,b)\in\Phi^+ \mbox{ and } p>0\}
\ee
We first prove that $V_+\neq \emptyset$.

Let $m\equiv [\frac{M}{2}]$.
The generators $S^{11}_{(1)},\ldots,S^{mm}_{(1)},S^{M+1,M+1}_{(1)},\ldots,
S^{M+n,M+n}_{(1)}$ form a Cartan subalgebra of $Osp(M|2n)$, so
there exists a least one eigenvector $v$ common to all $S^{aa}_{(1)}$ and
with eigenvalue $\mu=(\mu^{(1)}_1,\ldots,\mu^{(1)}_{M+2n})$.

If $v\in V_+$ then we have $V_{+}\neq\emptyset$.
If $v\not\in V_+$, by applying $S^{ab}_{(p)}$, $(a,b)\in\Phi^+$, to $v$ we 
obtain an other common eigenvector of the $S^{aa}_{(1)}$ with eigenvalue 
$\mu + \omega$, where $\omega$ is a $\ZZ_{>0}$-linear combination of 
the positive
roots. As $V$ is finite-dimensional, repeated applications of generators
$S^{ab}_{(n)}$, $(a,b)\in\Phi^+$, $n>0$, will lead to a non-vanishing vector
$v_+\in V$ such that
\be
S^{ab}_{(p)}v_+ = 0\;\forall (a,b)\in\Phi^+,\; p>0
\ee
\be
S^{aa}_{(1)}v_+ =\lambda^{(1)} v_+ \;\forall a
\ee
So $v_+\in V_+$ and $V_+$ contains at least one nonzero element.

One defines
\be
\ct_\pm=\{S^{ab}(u),\;\forall(a,b)\in\Phi^\pm \}
\ee
and $L$ (resp. $R$)  the left (resp. right) ideal 
generated by $\ct_{+}$ (resp. $\ct_{-}$). We also introduce the subalgebra
\be
\cy_{0}=\{y\in Y(M|2n)^+, \mbox{ such that } [S^{aa}_{(1)},y]=0\ 
\forall a=1,\ldots,M+2n\}
\ee
and correspondingly
\be
L_{0}=\cy_{0}\cap L \mb{and} R_{0}=\cy_{0}\cap R
\ee
Using the PBW theorem \ref{PBW}, one shows that $L_{0}=R_{0}\equiv 
I_{0}$ is a two-sided ideal so that $\cg=\cy_{0}/I_{0}$ is an algebra.
From the commutation relations (\ref{comSijSkl}), one gets that 
$[S^{aa}(u),S^{bb}(v)]\in I_{0}$, \ie $\cg$ is a commutative 
algebra.

By construction, $\forall v\in V_{+}$ and $i\in I_{0}$, one has 
$iv=0$, so that $V_{+}$ is a $\cg$-module. Since $\cg$ is 
commutative, there exists a nonzero common eigenvector 
$\xi\in V_{+}$. Now, 
let $V'=\cu(\ct_{-})\xi$: it is a non-zero submodule of $V$, and 
since $V$ is supposed irreducible, it must equal $V$. Thus, $\xi$ is 
a highest weight vector of $V$.

Finally, if there is another highest weight vector $\xi'$, the above 
construction ensures that $V=\cu(\ct_{-})\xi=\cu(\ct_{-})\xi'$ which 
is possible only for $\xi$ and $\xi'$ proportional.
\finprf

\begin{theor}[Necessary conditions for finite-dimensional irreps]\hfill\\
\label{theo-nec}
If the irreducible highest weight representation $V[\mu(u)]$ of  
$Y(M|2n)^+$ is finite-dimensional then the following relations hold:
\bea
&&\frac{\mu_{i}(u)}{\mu_{i+1}(u)}=\frac{P_{i+1}(u+1)}{P_{i+1}(u)} 
\mb{for} \left\{\begin{array}{l}
m+2\leq i\leq  M-1 \\[1.2ex]
M+n+1\leq i\leq  M+2n-1 
\end{array}\right.\\
&&\frac{\mu_{M+n+1}(-u)}{\mu_{M+n+1}(u)}=
\frac{P_{M+n+1}(u+1)P_{M+n+1}(-u)}{P_{M+n+1}(u)P_{M+n+1}(1-u)}
\eea

If $M=2m+1$, one among these two relations also holds: 
\be
\gamma(u)\frac{{\mu}_{m+1}(u)}{\mu_{m+2}(u)}=\frac{P_{m+1}(u+1)}{P_{m+1}(u)},
\mb{with}
\gamma(u)=\ 1\mb{or}
\frac{2u}{2u+1}
\label{eq.nec1}
\ee

If $M=2$, we have the supplementary condition
\be \displaystyle
\frac{{\mu}_{2}(-u)}{\mu_{2}(u)}=
\frac{P(u+1)P(-u)}{P(u)P(1-u)}\frac{(u+\gamma)(2u-1)}{u-\gamma)(2u+1)},
\mb{with}P(-\gamma)P(\gamma+1)\neq0
\label{eq.nec2}
\ee

Finally, for $M=2m>2$, we have the relations
\be\begin{array}{rl}
 \displaystyle
\frac{\mu^o_{m+1}(u)}{\mu_{m+2}(u)}&
\displaystyle
=\frac{P_{m+2}(u+1)}{P_{m+2}(u)}\\[1.2ex]
 \displaystyle
\gamma(u)\,\frac{{\mu^o}_{m+1}(-u)}{\mu^o_{m+1}(u)} &=
\displaystyle
\frac{P_{m+1}(u+1)P_{m+1}(-u)}{P_{m+1}(u)P_{m+1}(1-u)}
\end{array}
\label{eq.nec3}
\ee
with $\gamma(u)=1$ or $\frac{2u-1}{2u+1}$, and
$\mu^o_{m+1}(u)=\mu_{m+1}(u)$ or $\mu^\#_{m+1}(u)$.
We have introduce $\mu_{m+1}^\#(u)$, which is deduced from $\mu_{m+1}(u)$ 
by the 
action of the $\#$ automorphism (see definition \ref{defhash} 
and \cite{Molev} for more details).
\end{theor}
\prf
It is a direct consequence of the classification of finite-dimensional 
irreps for the algebras
$Y(M)^{\pm}$ done in \cite{Molev}. 
Since $Y(M)^{+}$ and $Y(2n)^-$ are subalgebras of $Y(M|2n)^+$, 
starting with an $Y(M|2n)^+$-irrep with highest 
weight $\xi$, and considering the cyclic span of $\xi$ with each of 
these subalgebras leads to the result.
\finprf

\subsection{Finite-dimensional irreps of $Y(1|2)^+$}
Let $V[\mu(u)]$ be an irrep 
of $Y(1|2)^+$ with 
highest weight 
$\mu(u)\equiv(\mu_1(u),\mu_2(u),\mu_3(u))$. From the symmetry relation 
(\ref{symmetry-mu}) we obtain that $\mu_1(u)$ is an even series in $u^{-1}$ 
and $\mu_{2}(u)$ can be deduced from $\mu_{3}(u)$.

\begin{prop}
\label{reecriture-mu1-mu3}
If $V[\mu_1(u),\mu_3(u)]$ is finite dimensional then there exists a formal 
even series 
$\psi(u)$ in $u^{-1}$ such that 
\bea
\label{mu1}
\psi(u)\mu_1(u) & = &(1-\alpha_1^2 u^{-2})\ldots(1-\alpha_k^2 u^{-2})\\
\label{mu3}
\psi(u)\mu_3(u) &=& (1-\alpha_1 u^{-1})\ldots(1-\alpha_k u^{-1})
(1+\beta_1 u^{-1})\ldots(1+\beta_k u^{-1})
\eea

\end{prop}
\prf:
Let $\xi$ be the highest weight vector of $V[\mu(u)]$.

Under $(S^{22}_{(1)},S^{33}_{(1)})$, $S^{32}_{(l_1)}\ldots S^{32}_{(l_s)}
S^{31}_{(p_1)}\ldots S^{31}_{(p_r)}\xi$ has weight 
$(\mu_2^{(1)}+2s+r,\mu_3^{(1)}-2s-r)$ whereas $S^{31}_{(i)}\xi$ 
has weight $(\mu_2^{(1)}+1,\mu_3^{(1)}-1)$. So $S^{31}_{(i)}\xi$ 
can only be written as a linear combination of vectors
$S^{31}_{(j)}\xi$. Let $k$ be the minimum 
non-negative integer such that $S^{31}_{(k+1)}\xi$ is a linear combination of 
the vectors $\xi_1\equiv S^{31}_{(1)}\xi,\ldots,\xi_k\equiv S^{31}_{(k)}\xi$ 
(such $k$ exists because $V[\mu(u)]$ is 
finite-dimensional).

We will prove that for any vector $S^{31}_{(r)}\xi$ with $r\geq k+1$ we have:
\be
\label{S31(r)}
S^{31}_{(r)}\xi=a_1^{(r)}\xi_1+\ldots+a_k^{(r)}\xi_k
\ee
where the $a_i^{(r)}$ are complex numbers. 

Equation (\ref{S31(r)}) 
is true for $r=k+1$ by definition of $k$.  
Taking $i=k=l=3$, $j=1$ in the commutation relation and exchanging 
$u$ and $v$, we get:

\bea
[S^{33}(u),S^{31}(v)] &=& -\frac{1}{u-v}(S^{31}(v)S^{33}(u)-
S^{31}(u)S^{33}(v))\\
&& -\frac{1}{u+v}S^{32}(u)S^{21}(v) +\frac{1}{u^2-v^2}S^{32}(u)S^{13}(v) 
\nonumber\\
&&+\frac{-1+u-v}{u^2-v^2}S^{32}(v)S^{13}(u)\nonumber
\eea
We multiply by $(u^2-v^2)$ and take the coefficient at $u^{0}v^{-p}$ 
($p\geq 1$). Using the fact that $S^{21}(u)\xi=S^{13}(u)\xi=0$ we obtain: 
\be
\label{recursion}
S^{33}_{(2)}S^{31}_{(p)}\xi=-S^{31}_{(p+1)}\xi+S^{31}_{(1)}S^{33}_{(p)}\xi+
S^{31}_{(p)}(S^{33}_{(2)}-S^{33}_{(1)})\xi 
\ee
For $i=1,\ldots,k-1$ (\ref{recursion}) gives:
\be
\label{S33(2)xii}
S^{33}_{(2)}\xi_i=-\xi_{i+1}+\mu_3^{(i)}\xi_1+(\mu_3^{(2)}-\mu_3^{(1)})\xi_i
\ee
For $i=k$, using $S^{31}_{(k+1)}\xi=a_1^{(k+1)}\xi_1+\ldots+
a_{k}^{(k+1)}\xi_k$ in (\ref{recursion}) gives:
\be
\label{S33(2)xik}
S^{33}_{(2)}\xi_k=-(a_1^{(k+1)}\xi_1+\ldots+a_{k}^{(k+1)}\xi_k)
+\mu_3^{(k)}\xi_1+(\mu_3^{(2)}-\mu_3^{(1)})\xi_k
\ee
So $\forall i\in\{1,\ldots,k\}$, $S^{33}_{(2)}\xi_i$ is a linear combination of
the $\{\xi_j\}_{j=1,\ldots,k}$.

Now suppose that $\forall r \in \{k+1,\ldots,p\}$ (where $p\geq k+1$), equation
(\ref{S31(r)}) holds. We then have:
\bea
S^{31}_{(p+1)}\xi &=& -S^{33}_{(2)}S^{31}_{(p)}+\mu_3^{(p)}\xi_1
   +(\mu_3^{(2)}-\mu_{3}^{(1)})S^{31}_{(p)}\xi \\
 &=& -\sum_{i=1}^{k}(a_i^{(p)}S^{33}_{(2)}\xi_i)
  +\mu_3^{(p)}\xi_1+(\mu_3^{(2)}-\mu_{3}^{(1)})
     \sum_{i=1}^{k}(a_i^{(p)}\xi_i)\nonumber
\eea
so $S^{31}_{(p+1)}\xi$ is a linear combination of the 
$\{\xi_j\}_{j=1,\ldots,k}$ and equation (\ref{S31(r)}) is proved 
by induction on $p$.
We can therefore write:
\be
\label{S31(u)}
S^{31}(u)\xi=a_1(u)\xi_1+\ldots+a_k(u)\xi_k
\ee
where $a_i(u)=u^{-i}+\sum_{s=k+1}^{\infty}a_i^{(s)}u^{-s}$.

We can rewrite (\ref{recursion}) as:
\be
\label{recursion2}
S^{33}_{(2)}S^{31}(v)\xi=-vS^{31}(v)\xi+\mu_3(v)\xi_1+
(\mu_3^{(2)}-\mu_3^{(1)})S^{31}(v)\xi 
\ee

On the other hand, applying $S^{33}_{(2)}$ on (\ref{S31(u)}) and using 
(\ref{S33(2)xii}) and (\ref{S33(2)xik}) we have:
\bea
\label{recursion3}
S^{33}_{(2)}S^{31}(v)\xi &=&\sum_{i=1}^{k}a_i(v)S^{33}_{(2)}\xi_i\\
&=& \left(\sum_{i=1}^{k}(a_i(v)\mu_3^{(i)})+
(\mu_3^{(2)}-\mu_3^{(1)})-a_k(v)a_1^{(k+1)}\right)\xi_1
\nonumber\\
& & +\sum_{i=2}^{k}\left(-a_{i-1}(v)+
(\mu_3^{(2)}-\mu_3^{(1)})a_i(v)-a_k(v)a_i^{(k+1)}\right)\xi_i\nonumber
\eea
Taking the coefficient at $\xi_i$ for $i=2,\ldots,k$ in (\ref{recursion2}) 
and (\ref{recursion3}) leads to:
\be
-a_{i-1}(v)+va_i(v)-a_i^{k+1}a_k(v)=0
\ee
so that:
\be
a_i(v)=a_k(v)\underbrace{(v^{k-i}-v^{k-1-i}a_k^{k+1}
                         -v^{k-i-2}a_{k-1}^{(k+1)}-\ldots -
                          va_{i+2}^{(k+1)}-a_{i+1}^{(k+1)}}_{A_i(v)}
\ee
for $i=1,\ldots,k-1$. The coefficient at $\xi_1$ in (\ref{recursion2}) 
and (\ref{recursion3}) leads to:
\be
\mu_3(v)=\underbrace{\sum_{i=1}^{k}\mu_3^{(i)}a_k(v)A_i(v)}_{
\mbox{pol. of degree $k-1$} } -a_1^{(k+1)}a_k(v)+
\underbrace{va_k(v)A_{1}(v)}_{\mbox{monic pol. of degree $k$} } 
\ee
So $\mu_3(v)=a_k(v)B(v)$ where B(v) is a monic polynomial in $v$ of 
degree $k$.

For $\mu_1(v)$ we use:
\bea
[S^{11}(u),S^{31}(v)] & =& \frac{u+v+1}{u^2-v^2} S^{31}(u)S^{11}(v)
 - \frac{1}{u+v} S^{12}(u)S^{11}(v)\nonumber\\
 && - \frac{2v+1}{u^2-v^2} S^{31}(v)S^{11}(u)\nonumber
\eea
Notice that $S^{11}(u)$ is an even series in $u^{-1}$. Using the same procedure
we find that $\mu_1(v)=a_k(v)C(v)$ where C(v) is a monic polynomial 
in $v$ of degree $k$.

Defining $\varphi(u)=(a_k(u)u^k)^{-1}$, we have:
\bea
\varphi(u)\mu_1(u) &=&(1+\alpha_1 u^{-1})\ldots(1+\alpha_k u^{-1})\\
\varphi(u)\mu_3(u) &=&(1+\beta_1 u^{-1})\ldots(1+\beta_k u^{-1})
\eea
where the $\alpha_i$'s and the $\beta_i$'s are complex numbers. 

The formal series
\bea
\psi(u)&\equiv&\varphi(u)(1-\alpha_1 u^{-1})\ldots(1-\alpha_k u^{-1})\\
&=&\frac{(1-\alpha_1^2 u^{-2})\ldots(1-\alpha_k^2 u^{-2})}{\mu_1(u)}
\nonumber
\eea
is an even series in $u^{-1}$. The composition of the automorphism 
$S(u)\rightarrow \psi_(u)S(u)$ with $V[\mu(u)]$ is an irrep 
 with the following highest weight which we shall 
again denote by $\mu(u)$
\bea
\mu_1(u) & = &(1-\alpha_1^2 u^{-2})\ldots(1-\alpha_k^2 u^{-2})\\
\mu_3(u) &=& (1-\alpha_1 u^{-1})\ldots(1-\alpha_k u^{-1})
(1+\beta_1 u^{-1})\ldots(1+\beta_k u^{-1})
\eea
\finprf
Thus, up to an automorphism of $Y(1|2)^+$, we can assume that $\mu_1(u)$ 
and $\mu_3(u)$ are polynomials in $u^{-1}$.
\begin{theor}\label{theoY12}
Let $V[\mu_1(u),\mu_3(u)]$ be an irrep of $Y(1|2)^+$. Suppose $\mu_1(u)$
and $\mu_3(u)$ satisfy
\bea
\label{mu1-mu3}
\frac{{\mu}_1(u)}{\mu_3(u)}&=&\frac{P(u+1)}{P(u)} \frac{R(u)}{Q(u)}\\
\label{mu3-mu3}
\frac{\mu_3(-u)}{\mu_3(u)}&=&\frac{P(u+1)P(-u)}{P(u)P(1-u)}
\eea
where $P(u)$, $Q(u)$ and $R(u)$ are a monic polynomial, $Q(u)$ and $R(u)$
are even in $u$ and of same degree. 

Then $V$ is finite-dimensional.
\end{theor}
\prf
We call $p$ (resp. $2r$) the degree of $P(u)$ (resp. $Q(u)$ and
 $R(u)$). Since $R(u)$ and $Q(u)$ are even, they write
\be
R(u)=R_0(u)R_0(-u)\mb{;} Q(u)=Q_0(u)Q_0(-u)\mb{with}dg(R_0)=dg(Q_0)=r
\label{polyRQ}
\ee

We introduce:
\bea
\lambda_1(u)&=&u^{-s-r}P(u+1)R_0(u)\\
\lambda_2(u)&=&u^{-s-r}P(u+1)Q_0(u)\\
\lambda_3(u)&=&u^{-s-r}P(u)Q_0(u)
\eea
Let $L[\lambda(u)]$ be the corresponding irreducible highest weight 
module of $Y(1|2)$. 
Since $\lambda_1(u)/\lambda_2(u)=R_0(u)/Q_0(u)$ and 
$\lambda_2(u)/\lambda_3(u)=P(u+1)/P(u)$, according to \cite{zhang} 
$L[\lambda(u)]$ is finite-dimensional. The cyclic $Y(1|2)^+$-span of its 
highest weight 
vector is a finite-dimensional representation $V[\mu'(u)]$ of $Y(1|2)^+$ 
with $\mu'_1(u)=\lambda_1(u)\lambda_1(-u)$ and 
$\mu'_3(u)=\lambda_3(u)\lambda_2(-u)$. By construction, 
the polynomials  $\mu'_i(u)$ 
satisfy  (\ref{mu1-mu3}-\ref{mu3-mu3}). This implies that:
\be
\psi(u)\equiv\frac{\mu_3(u)}{\mu'_3(u)}=\frac{\mu_3(-u)}{\mu'_3(-u)}
\ee
is an even series in $u^{-1}$ and
\be
\mu_1(u)=\frac{\mu_3(u)}{\mu'_3(u)}\mu'_1(u)=\psi(u)\mu'_1(u)
\ee
Thus, there exists an automorphism $S(u)\rightarrow\psi(u)S(u)$ of $Y(1|2)^+$
such that its composition with the representation $V[\mu'(u)]$ is isomorphic 
to $V[\mu(u)]$: $V[\mu(u)]$ is therefore finite-dimensional.
\finprf
\begin{guess}\label{guess1}
The sufficient condition (\ref{mu1-mu3}) of theorem \ref{theoY12} for the 
existence of 
finite-dimensional irreps of $Y(1|2)^+$ is also a 
necessary condition.
\end{guess}
We remind that the condition (\ref{mu3-mu3}) has been proved to be necessary 
(see theorem \ref{theoY12}), so that the conjecture \ref{guess1} just says 
that the theorem \ref{theoY12} states necessary and sufficient conditions 
for finite-dimensional irreps of $Y(1|2)^+$.

\subsection{The general case $Y(2m+1|2n)^+$}
\begin{theor}
\label{theodd}
Let $V=V[\mu_{m+1}(u),...\mu_{2m+1}(u),\mu_{M+n+1}(u),...,\mu_{M+2n}(u)]$ be 
an irrep of $Y(2m+1|2n)^+$. We take $m\geq1$ and note $M=2m+1$.

Suppose the weights $\mu_i(u)$ obey
\bea
&&\frac{\mu_{i}(u)}{\mu_{i+1}(u)}=\frac{P_{i+1}(u+1)}{P_{i+1}(u)} 
\mb{for} \left\{\begin{array}{l}
 m+2\leq i\leq  2m \\[1.2ex]
M+n+1\leq i\leq  M+2n-1 
\end{array}\right.
\label{cond1}\\
&&\frac{\mu_{M+n+1}(-u)}{\mu_{M+n+1}(u)}=
\frac{P_{M+n+1}(u+1)P_{M+n+1}(-u)}{P_{M+n+1}(u)P_{M+n+1}(1-u)}
\label{cond2}\\
&&\gamma(u)\frac{{\mu}_{m+1}(u)}{\mu_{m+2}(u)}=\frac{P_{m+2}(u+1)}{P_{m+2}(u)}
\label{cond3}\\
&&\frac{{\mu}_{m+1}(u)}{\mu_{M+n+1}(u)}=
\frac{P_{M+n+1}(u+1)}{P_{M+n+1}(u)}\frac{R(u)}{Q(u)}\label{cond4}
\eea
{with}
$R(u)$ and $Q(u)$ even and of same degree.
In the above formulas, $\gamma(u)$ 
is either $1$, and the corresponding relations 
 will be called case (a), or 
$\frac{2u}{2u+1}$, case (b).

Then $V$ is finite-dimensional.

\null

Under the assumption of conjecture \ref{guess1}, the above 
sufficient conditions
are also necessary ones.
\end{theor}
\prf
First, let case (a) hold. We note $s_i$ the degree of
the polynomials ${P_i(u)}$, decompose $R(u)$ and $Q(u)$ as in 
(\ref{polyRQ}).
and note:
\be
P_+(u) = \prod_{a=m+2}^{M}\!P_a(u)\mb{;}
P_-(u) = \prod_{a=M+n+1}^{M+2n}\!P_a(u)\mb{;}
s_0=r+\!\sum_{i=m+1}^{M}\!s_i+\!\sum_{i=M+n+1}^{M+2n}\!\!s_i
\ee
We also define:
\be
\lambda_i(u) = u^{-s_0}P_+(u+1)P_-(u+1)R_0(u),\ \ \ i=1,\ldots,m+1
\label{L1}
\ee
\be
\lambda_i(u)=u^{-s_0}P_-(u+1)R_0(u)
\prod_{a=m+2}^{i}\! P_a(u)\prod_{a=i+1}^{2m+1}
\!P_a(u+1),\ \ \ i=m+2,\ldots,M 
\label{L2}
\ee
\be
\lambda_i(u) = u^{-s_0}P_+(u+1)P_-(u+1)Q_0(u),\ \ \ i=M+1,\ldots,M+n
\label{L3}
\ee
\bea
\lambda_i(u)& =& u^{-s_0}P_+(u+1)Q_0(u)\prod_{a=M+n+1}^{i} 
\!\!P_a(u)\prod_{a=i+1}^{M+2n}\!P_a(u+1),\nonu
&& i=M+n+1,\ldots,M+2n
\label{L4}
\eea

We therefore have:
\bea
\label{lambda1}
\frac{\lambda_i(u)}{\lambda_{i+1}(u)} &=&\frac{P_{i+1}(u+1)}{P_{i+1}(u)}
\;\mbox{ for }\;\left\{\begin{array}{l} 
i=m+1,\ldots,M-1\\ i=M+n,\ldots,M+2n-1
\end{array}\right.
\\
\label{lambda2}
\frac{\lambda_i(u)}{\lambda_{i+1}(u)}&=&1
\;\mbox{ for }\; \left\{\begin{array}{l} 
i=1,\ldots,m\\ i=M+1,\ldots,M+n-1
\end{array}\right.
\\
\label{lambda3}
\frac{\lambda_{M}(u)}{\lambda_{M+1}(u)}&=&\frac{R_0(u)}{Q_0(u)}
\,\frac{P_+(u)}{P_+(u+1)}
\eea
We consider the  highest weight irrep $L[\lambda(u)]$ 
of $Y(2m+1|2n)$.
According to property \ref{Irreps-YMN}, the relations 
(\ref{lambda1}), (\ref{lambda2}) and (\ref{lambda3}) 
ensure that $L[\lambda(u)]$ is finite-dimensional. 
The cyclic $Y(2m+1|2n)^+$-span of its highest weight vector is a 
finite-dimensional representation
with highest weights $\mu'_i(u)=
\lambda_i(u)\lambda_{\bar{\imath}}(-u)$ for $i=m+1,\ldots,2m+1$ and 
$i=M+n+1,\ldots,M+2n$. Its irreducible quotient is a finite-dimensional
irrep $V[\bar\mu(u)]$ of $Y(2m|2n)^+$.

Moreover, the $\mu'_i(u)$ verify:
\bea
\frac{\mu'_i(u)}{\mu'_{i+1}(u)} &=& \frac{P_{i+1}(u+1)}{P_{i+1}(u)}=
\frac{\mu_i(u)}{\mu_{i+1}(u)}\ ,\ \ 
 \left\{\begin{array}{l}
i=m+1,\ldots,M-1 \\ i=M+n+1,\ldots,M+2n-1 
\end{array}\right.\ \ \ \ \ \\
\frac{\mu'_{M+n+1}(-u)}{\mu'_{M+n+1}(u)}&=&\frac{P_{M+n+1}(u+1)P_{M+n+1}(-u)}
{P_{M+n+1}(u)P_{M+n+1}(1-u)}=\frac{\mu_{M+n+1}(-u)}{\mu_{M+n+1}(u)}
\\
\frac{\mu'_{m+1}(u)}{\mu'_{M+n+1}(u)}&=&\frac{P_{M+n+1}(u+1)}{P_{M+n+1}(u)}
\frac{R(u)}{Q(u)}
=\frac{\mu_{m+1}(u)}{\mu_{M+n+1}(u)}
\eea
The formal series
\be
\psi(u)=\frac{\mu_{M+n+1}(-u)}{\mu'_{M+n+1}(u)}
\ee
is an even series in $u^{-1}$ and we have
\be
\frac{\mu_i(u)}{\mu'_i(u)}=\psi(u),\ \ \forall i
\ee
Hence there exists an automorphism $S(u)\rightarrow \psi(u)S(u)$ 
of $Y(2m+1|2n)^+$ such that its composition with $V[\mu'(u)]$ is 
isomorphic to $V[\mu(u)]$. This later is thus finite-dimensional.

\null 

Now, let case (b) hold.
We introduce the $osp(2m+1|2n)$ representation $V_0$,
of highest weight $l_i=-\half$ for $i=m+2,...,2m+1$ and $l_i=0$ for 
$i=m+1,M+n+1,...,M+2n$
and promote it to a $Y(2m+1|2n)^+$ 
representation using the evaluation map. The corresponding highest weight 
 has components $l_i(u)=\frac{2u}{2u+1}$  for $i=m+2,...,2m+1$,
and $l_{i}(u)=1$ for 
$i=m+1,M+n+1,...,M+2n$.

Moreover, making the same construction as for case (a), 
and considering the tensor product 
$L[\lambda(u)]\otimes V_0$, we get a finite dimensional representation 
$V[\mu''(u)]$ obeying the relations of case (b). Its irreducible subquotient 
is isomorphic to $V[\mu(u)]$, which is therefore finite-dimensional.

\null

Conversely, let us suppose that the irrep $V[\mu(u)]$ is finite dimensional. 
From theorem \ref{theo-nec}, one already knows that the conditions 
(\ref{cond1}), (\ref{cond2}) and (\ref{cond3}) must be satisfied.

Suppose also that the conjecture \ref{guess1} holds.
The subalgebra generated by the coefficients of $S^{ij}(u)$, 
$i,j=m+1,M+n,M+n+1$ is isomorphic to $Y(1|2)^+$. The cyclic span
of the highest weight vector of $V[\mu(u)]$ with respect to this subalgebra
is a representation with highest weight 
$(\mu_{m+1}(u),\mu_{M+n+1}(u))$.
Its irreducible quotient is finite-dimensional and so, we have relation 
(\ref{cond4}). 
\finprf

\subsection{Finite-dimensional irreps of $Y(2|2)^+$}
Let $V=V[\mu(u)]$ be an irrep of $Y(2|2)^+$ with 
highest weight 
$\mu(u)$.

\begin{prop}
\label{reecriture-mu2-mu4}
If $V[\mu(u)]$ is finite dimensional then there exists a formal even series 
$\psi(u)$ in $u^{-1}$ such that 
\bea
\label{mu2}
\psi(u)\mu_2(u) &=&
(1-\alpha_1 u^{-1})\ldots(1-\alpha_k u^{-1})\\
\label{mu4}
\psi(u)\mu_4(u) &=& (1-\beta_1 u^{-1})\ldots(1-\beta_k u^{-1})
\eea
\end{prop}
\prf
The proof is very similar to the case of $Y(1|2)^+$, and we leave it to the 
reader.
Note that the calculation being achieved using the fermionic generator 
$S_{13}(u)$ (instead of the even bosonic one $S_{12}(u)$), there is
no difference in the proof for $Y(1|2)^+$ and $Y(2|2)^+$, in
opposition with the $Y(2)^+$ and $Y(2)^-$ cases \cite{Molev}.
\finprf
\begin{theor}\label{theoY22}
Let $V[\mu_2(u),\mu_4(u)]$ be an irrep of $Y(2|2)^+$. 
If $\mu_2(u)$
and $\mu_4(u)$ satisfy
\bea
\label{mu2-mu4}
\frac{\mu_2(u)}{\mu_4(u)}&=&\frac{u-\gamma}{u+\half}\,
\frac{P_4(u+1)\,P_2(u)}{P_4(u)\,P_2(u+1)}\,\frac{R(u)}{Q(u)}
 \\
\label{mu4-mu4}
\frac{\mu_4(-u)}{\mu_4(u)}&=& \left(\frac{u+\half}{u-\half}\right)^2\,
\frac{P_4(u+1)\,P_4(-u)}{P_4(u)\,P_4(1-u)}
\eea
then $V$ is finite dimensional.

In the above formulas,
 $P_2(u)$, $P_4(u)$, $Q(u)$ and $R(u)$ are monic polynomials,
 $Q(u)$ and $R(u)$ are even in $u$ and of same degree, and $\gamma\in\CC$. 
\end{theor}
\prf
We call $p_2$ (resp. $p_4$, resp. $2r$) the degree of $P_2(u)$ 
(resp. $P_4(u)$, resp. $Q(u)$ and
 $R(u)$) and decompose $Q(u)$ and $R(u)$ as in (\ref{polyRQ}).
Let $L[\lambda(u)]$ be the irrep of $Y(2|2)$ with
weights
\bea
\lambda_1(u) &=& u^{-p_2-p_4-r}P_4(u+1)P_2(u+1)R_0(u)\\
\lambda_2(u) &=& u^{-p_2-p_4-r}P_4(u+1)P_2(u)R_0(u)\\
\lambda_3(u) &=& u^{-p_2-p_4-r}P_4(u+1)P_2(u+1)Q_0(u)\\
\lambda_4(u) &=& u^{-p_2-p_4-r}P_4(u)P_2(u+1)Q_0(u)
\eea
Since we have 
\be
\frac{\lambda_1(u)}{\lambda_2(u)}=\frac{P_2(u+1)}{P_2(u)}\mb{;}
\frac{\lambda_2(u)}{\lambda_3(u)}=\frac{P_2(u+1)R_0(u)}{P_2(u)\,Q_0(u)}\mb{;}
\frac{\lambda_3(u)}{\lambda_4(u)}=\frac{P_4(u+1)}{P_4(u)}
\ee
$L[\lambda(u)]$ is finite dimensional. 
The $Y(2|2)^+$-cyclic span of its 
highest weight vector is a finite dimensional $Y(2|2)^+$-representation 
$V[\mu'(u)]$ of weights $\mu_2'(u)=\lambda_2(u)\lambda_1(-u)$ and 
$\mu_4'(u)=\lambda_4(u)\lambda_3(-u)$. These weights obey the relations 
\bea
\frac{\mu'_2(u)}{\mu'_4(u)}&=&
\frac{P_4(u+1)\,P_2(u)}{P_4(u)\,P_2(u+1)}\,\frac{R(u)}{Q(u)}
 \\
\frac{\mu'_4(-u)}{\mu'_4(u)}&=& 
\frac{P_4(u+1)\,P_4(-u)}{P_4(u)\,P_4(1-u)}
\eea

We now consider the $osp(2|2)$ finite-dimensional irrep $V_0$ with weights
$l_2=-\gamma-\half$ and $l_4=-1$. Through the evaluation map, its
provides a finite-dimensional representation of $Y(2|2)^+$ with weights 
\be
l_2(u)=\frac{u-\gamma}{u+\half}\mb{and}
l_4(u)=\frac{u-\half}{u+\half}
\ee
The tensor product $L[\lambda(u)]\otimes V_0$ is thus a finite-dimensional 
representation of $Y(2|2)^+$, with weights $\mu_i''(u)=\mu_i'(u)l_i(u)$,
$i=2,4$. They obey to relations
(\ref{mu2-mu4}) and (\ref{mu4-mu4}), so that the irreducible quotient 
provide a finite-dimensional irrep isomorphic to $V$. 
Thus, $V$ finite-dimensional.
\finprf

Note that the polynomial $P(u)=(u-\half)^2$ obeys the relation 
$P(1-u)=P(u)$, so that the condition on $\mu_4(u)$ does not differ from 
the one obtained for $Y(2)^-$.

\begin{guess}\label{guess2}
The sufficient condition (\ref{mu2-mu4}) of theorem \ref{theoY22} for the 
existence of 
finite-dimensional irreps of $Y(2|2)^+$ is also a 
necessary condition.
\end{guess}

\subsection{The general case $Y(2m|2n)^+$}

\begin{theor}[Case of $Y(2|2m)^+$]\hfill\\
Let $V=V[\mu_{2}(u),\mu_{n+3}(u),...,\mu_{2+2n}(u)]$ be 
an irrep of $Y(2|2n)^+$. 

Suppose the weights $\mu_i(u)$ obey
\bea
&&\frac{\mu_{i}(u)}{\mu_{i+1}(u)}=\frac{P_{i+1}(u+1)}{P_{i+1}(u)} 
\mb{for}
n+3\leq i\leq  2n+1 
\label{dit1}\\
&&\frac{\mu_{n+3}(-u)}{\mu_{n+3}(u)}=
\left(\frac{u-\half}{u+\half}\right)^2\,
\frac{P_{n+3}(u+1)P_{n+3}(-u)}{P_{n+3}(u)P_{n+3}(1-u)}
\label{dit2}\\
&&\frac{{\mu}_{2}(u)}{\mu_{n+3}(u)}=\frac{u-\gamma}{u+\half}\,
\frac{P_{n+3}(u+1)\,P_{2}(u)}{P_{n+3}(u)\,P_{2}(u+1)}
\,\frac{R(u)}{Q(u)}\label{dit3}
\eea
{with}
$R(u)$ and $Q(u)$ even and of same degree, and $\gamma\in\CC$.

Then $V$ is finite-dimensional.

\null

Under the assumption of conjecture \ref{guess2}, the conditions 
(\ref{dit1})-(\ref{dit3}) are
necessary and sufficient conditions for $V$ to be a finite-dimensional irrep.
\end{theor}
\prf
The proof is similar to the previous ones. One constructs a 
finite-dimensional irrep for $Y(2|2n)^+$ which fulfils the 
conditions (\ref{dit1})-(\ref{dit3}). It takes the form 
$V'=L[\lambda(u)]\otimes V_\gamma$. $L[\lambda(u)]$ is constructed as in
 theorem \ref{theodd}. The  
 finite-dimensional $osp(2|2n)$-irrep $V_\gamma$ has weight
\be
l_2=-\gamma-\half\mb{and} l_i=-1, \mb{for} n+3\leq i\leq 2n
\ee
Looking at the $osp(2|2n)$-span of the highest weight vector
and taking the irreducible subquotient, we get a finite-dimensional irrep
$V'$ with highest weight obeying (\ref{dit1})-(\ref{dit3}). $V$ being 
isomorphic to $V'$, it is therefore finite dimensional.

Conversely, assuming the conjecture \ref{guess2}, and looking at 
the subalgebras
$Y(2|2)^+$ and $Y(2n)^-$, one easily proves that the conditions
(\ref{dit1})-(\ref{dit3}) are necessary conditions.
\finprf

\begin{theor}[Case of $m>1$]\hfill\\
\label{theven}
Let $V=V[\mu_{m+1}(u),...\mu_{M}(u),\mu_{M+n+1}(u),...,\mu_{M+2n}(u)]$ be 
an irrep of $Y(2m|2n)^+$. We note $M=2m$ and take $m>1$.

Suppose the weights $\mu_i(u)$ obey
\bea
&&\frac{\mu_{i}(u)}{\mu_{i+1}(u)}=\frac{P_{i+1}(u+1)}{P_{i+1}(u)} 
\mb{for} \left\{\begin{array}{l}
 m+1\leq i\leq  2m-1 \\[1.2ex]
M+n+1\leq i\leq  M+2n-1 
\end{array}\right.
\label{tion1}\\
&&\frac{\mu_{M+n+1}(-u)}{\mu_{M+n+1}(u)}=
\frac{P_{M+n+1}(u+1)P_{M+n+1}(-u)}{P_{M+n+1}(u)P_{M+n+1}(1-u)}
\label{tion2}\\
&&\gamma(u)\frac{{\mu}_{m+1}(u)}{\mu_{M+n+1}(u)}=
\frac{P_{M+n+1}(u+1)\,P_{m+1}(u)}{P_{M+n+1}(u)\,P_{m+1}(u+1)}
\,\frac{R(u)}{Q(u)}\label{tion3}
\eea
{with}
$R(u)$ and $Q(u)$ even and of same degree, and 
$\gamma(u)=1$ or 
$\gamma(u)=\frac{2u-1}{2u+1}$.

Then $V$ is finite-dimensional.
\end{theor}
\prf
We start with the case $\gamma(u)=1$, and do the
same construction as in theorem \ref{theodd}, to get weights 
$\lambda_i(u)$ defined as in equations (\ref{L1})-(\ref{L4}),
with now $M=2m$. We introduce:
\bea
\lambda_i'(u) &=& P_{m+1}(u+1)\lambda_i(u)\mb{for} \left\{
\begin{array}{l} i=1,...,m\\ i=M+1,...,M+2n \end{array}
\right.\\
\lambda_i'(u) &=& P_{m+1}(u)\lambda_i(u)\mb{for} 
 i=m+1,...,M
\eea
For these new weights, the relations (\ref{lambda1})-(\ref{lambda3}) 
are still valid when 
$i\neq m,M$. In these later cases, we get
\bea
\frac{\lambda'_m(u)}{\lambda'_{m+1}(u)} &=& \frac{P_{m+1}(u+1)}{P_{m+1}(u)}\\
\frac{\lambda'_M(u)}{\lambda'_{M+1}(u)} &=& \frac{P_{m+1}(u)}{P_{m+1}(u+1)}
\frac{P_{+}(u)}{P_{+}(u+1)}\frac{R_0(u)}{Q_0(u)}
\eea
Thus, the $Y(2m|2n)$-irrep $L[\lambda'(u)]$ is still finite-dimensional.
The cyclic $Y(2m|2n)^+$-span of its highest weight vector is a representation
with highest weight $\bar\mu_i(u)=
\lambda'_i(u)\lambda'_{\bar{\imath}}(-u)$ for $i=m+1,\ldots,2m$ and 
$i=M+n+1,\ldots,M+2n$. Its irreducible quotient is a finite-dimensional
irrep $V[\bar\mu(u)]$ of $Y(2m|2n)^+$.

Moreover, the weights $\bar\mu_i(u)$ 
for $i=m+1,...,M-1$ on the one hand, and $i=M+n+1,..,M+2n-1$ on the other hand,
verify the same relations
as the $\mu'_i(u)$, i.e. the conditions (\ref{tion1}-\ref{tion2}).
For the remaining relation, one computes;
\be
\frac{\bar\mu_{m+1}(u)}{\bar\mu_{M+n+1}(u)} = 
\frac{P_{m+1}(u)}{P_{m+1}(u+1)}\,\frac{\mu'_{m+1}(u)}{\mu'_{M+n+1}(u)}
\ee
which gives the relation (\ref{tion3}). 

The weights $\mu_i(u)$ and $\bar\mu_i(u)$ obeying both the
relations (\ref{tion1})-(\ref{tion3}),
there exists an automorphism $S(u)\rightarrow \psi(u)S(u)$ 
of $Y(2m|2n)^+$ such that its composition with $V[\bar\mu(u)]$ is 
isomorphic to $V[\mu(u)]$. This later is thus finite-dimensional.

\null 

If now $\gamma(u)=\frac{2u}{2u-1}$, we construct 
the tensor product of the above representation by the $osp(2m|2n)$ 
finite-dimensional irrep $V_0$ with weights $l_i=-\half$ for 
$m+1\leq i\leq 2m$ and $l_i=-1$ for $2m+n+1\leq i\leq 2m+2n$. $V_0$
provides a finite-dimensional representation for $Y(2m|2n)^+$ with weights
$l_i(u)=\frac{2u}{2u+1}$ for  $m+1\leq i\leq 2m$ and $l_i(u)=\frac{2u}{2u-1}$
 for $2m+n+1\leq i\leq 2m+2n$. The weights of the tensor product obey
the relations (\ref{tion1})-(\ref{tion3}), and we
conclude as in theorem \ref{theodd}.
\finprf 

{\bf Remark:}

Conversely, let us suppose that the irrep $V[\mu(u)]$ is finite dimensional
and that the conjecture \ref{guess2} holds. 
From theorem \ref{theo-nec}, one already knows that the conditions 
(\ref{tion1}) and (\ref{tion2}) must be satisfied (for $i\neq m+1$). 
Moreover, we get also 
\be
\gamma(u)\,\frac{{\mu^o}_{m+1}(-u)}{\mu^o_{m+1}(u)}=
\frac{P_{m}(u+1)P_{m}(-u)}{P_{m}(u)P_{m}(1-u)}
\ee
with $\gamma(u)$ and $\mu^o_{m+1}(u)$
defined as in the theorem \ref{theo-nec}. 

We suppose also that $\mu^o_{m+1}(u)=\mu_{m+1}(u)$, which turns 
out to suppose that (\ref{tion1}) is valid for $i=m+1$.

The subalgebra generated by the coefficients of $S^{ij}(u)$, 
$i,j=m,m+1,M+n,M+n+1$ is isomorphic to $Y(2|2)^+$. The cyclic span
of the highest weight vector of $V[\mu(u)]$ with respect to this subalgebra
is a representation with highest weight 
$(\mu_{m+1}(u),\mu_{M+n+1}(u))$.
Its irreducible quotient is finite-dimensional and so, we have 
\be
\frac{\mu_{m+1}(u)}{\mu_{M+n+1}(u)}=\frac{u-\gamma}{u+\half}\,
\frac{P_{M+n+1}(u+1)\,P_{m+1}(u)}{P_{M+n+1}(u)\,P_{m+1}(u+1)}\,
\frac{R(u)}{Q(u)}
\ee
We look at the $osp(2m|2n)$ irrep induced by the generators 
$S^{ab}_{(1)}$ acting on the highest weight vector. It is finite dimensional,
so that we must have
\bea
&&l_{i+1}-l_i\in\ZZ_+\mb{for} \left\{\begin{array}{l} m+2\leq i\leq 2m\\
M+n+1\leq i\leq M+2n-1
\end{array}\right.\\
&&-(l_{m+2}+l_{m+1})\in\half\ZZ_+\mb{and} -l_{M+n+1}\in\ZZ_+\\
&&l_{m+1}-l_{M+n+1}\in\half\ZZ_+
\eea
where $\mu_i(u)=1+u^{-1}l_i+...$.
The above relations (on the weights $\mu_i(u)$)
imply the following constraints:
\bea
&&l_{i+1}-l_i\in\ZZ_+\mb{for} \left\{\begin{array}{l} m+2\leq i\leq 2m\\
M+n+1\leq i\leq M+2n-1
\end{array}\right.\\
&&-(l_{m+1}+l_{m+2})\in\ZZ_+\mb{and} -l_{M+n+1}\in\ZZ_+\\
&&l_{m+1}-l_{M+n+1}\in\ZZ_++(-\gamma-\half)
\eea
This implies in particular that $\gamma\in-\half\ZZ_+$,
so that we are back to the conditions (\ref{tion3}).
Thus, we are led to the following:
\begin{guess}
The sufficient conditions of theorem \ref{theven} for the 
existence of 
finite-dimensional irreps of $Y(2m|2n)^+$ are also 
necessary conditions.
\end{guess}

\sect{Conclusion\label{conc}}
We have defined the notion of twisted Yangians in the context of 
superalgebra $gl(M|N)$. It appears that most of the properties of the
twisted Yangians $Y^\pm(N)$ can be exhibited in the superalgebra case.
However, only $Y^+(M|2n)$ and $Y^-(2m|N)$ can be defined, and appear to be 
isomorphic. We thus concentrate on $Y^+(M|2n)$. Its
 finite dimensional irreducible representations
have being studied. 
$Y^+(M|2n)$ is also a coideal subalgebra of $Y(M|2n)$,
and is a deformation of the polynomial superalgebra $\cu(osp(M|2n)[x])$.

From a mathematical point of view, the centre of this algebra remains 
to be studied, and in particular the
notion of Sklyanin determinant (which appear in the context of twisted 
Yangians) has to be generalised to this case. Note that the notion of 
quantum Berezinian, which generates central element of $Y(M|N)$ has 
been introduced in \cite{naz}.

From the physical of view, the integrable systems with
boundary that could be relevant for such an algebra has to be determined.
Nonlinear sigma models based on a supergroup seem to be a good candidate.

\null

{\bf Acknowledgements}

We warmly thank D. Arnaudon and L. Frappat for fruitful
 comments.

A part of the redaction was done at the School of Mathematics, 
Sydney university, and E.R. thanks them for hospitality.

\end{document}